\documentclass[12pt]{amsart}

\usepackage{fix-cm}
\usepackage[normalem]{ulem}

\usepackage{lineno}

\usepackage{soul}

\usepackage{graphicx}

\usepackage{amssymb, amsmath, amsthm, hyperref, euscript, color}

\usepackage[dvipsnames]{xcolor}

\usepackage{amscd}

\usepackage{tikz-cd}

\usepackage{bbm}

\usepackage{enumitem}

\usepackage[english]{babel}

\usepackage{mathtools}

\usepackage{makecell}

\numberwithin{equation}{section}

\newtheoremstyle{theor}{6pt plus 1pt minus 1pt}{6pt plus 1pt minus 1pt}{\slshape}{}{\bfseries}{.}{5pt plus 1pt minus 1pt}{}
\newtheoremstyle{def}{6pt plus 1pt minus 1pt}{6pt plus 1pt minus 1pt}{}{}{\bfseries}{.}{5pt plus 1pt minus 1pt}{}
\newtheoremstyle{rmk}{6pt plus 1pt minus 1pt}{6pt plus 1pt minus 1pt}{}{}{\bfseries}{.}{5pt plus 1pt minus 1pt}{}
\newtheoremstyle{claim}{6pt plus 1pt minus 1pt}{6pt plus 1pt minus 1pt}{}{}{\bfseries}{.}{5pt plus 1pt minus 1pt}{}

\theoremstyle{theor}
\newtheorem{newstatement}{newstatement}
\newtheorem{lemma}[newstatement]{Lemma}
\newtheorem{theorem}[newstatement]{Theorem}
\newtheorem*{theorem*}{Theorem 2}
\newtheorem{corollary}[newstatement]{Corollary}
\newtheorem{proposition}[newstatement]{Proposition}

\theoremstyle{def}
\newtheorem{definition}[newstatement]{Definition}

\theoremstyle{rmk}
\newtheorem{remark}[newstatement]{Remark}
\newtheorem{example}[newstatement]{Example}
\newtheorem*{example*}{Example}

\theoremstyle{claim}

\theoremstyle{theor}
\newtheorem{thm}{Theorem}

\expandafter\let\expandafter\oldproof\csname\string\proof\endcsname
\let\oldendproof\endproof
\renewenvironment{proof}[1][\proofname]{%
  \oldproof[\slshape #1]%
}{\oldendproof}
\def\provedboxcontents#1{$\square$}

\makeatletter
\newsavebox\myboxA
\newsavebox\myboxB
\newlength\mylenA

\newcommand*\xoverline[2][0.75]{%
    \sbox{\myboxA}{$\m@th#2$}%
    \setbox\myboxB\null
    \ht\myboxB=\ht\myboxA%
    \dp\myboxB=\dp\myboxA%
    \wd\myboxB=#1\wd\myboxA
    \sbox\myboxB{$\m@th\overline{\copy\myboxB}$}
    \setlength\mylenA{\the\wd\myboxA}
    \addtolength\mylenA{-\the\wd\myboxB}%
    \ifdim\wd\myboxB<\wd\myboxA%
       \rlap{\hskip 0.5\mylenA\usebox\myboxB}{\usebox\myboxA}%
    \else
        \hskip -0.5\mylenA\rlap{\usebox\myboxA}{\hskip 0.5\mylenA\usebox\myboxB}%
    \fi}
\makeatother

\newcommand{\Q}{\mathbb{Q}}

\newcommand{\R}{\mathbb{R}}

\newcommand{\Z}{\mathbb{Z}}

\newcommand{\N}{\mathbb{N}}

\newcommand{\C}{\mathbb{C}}

\newcommand{\CP}{{\mathbb C\mkern-0.5mu\mathrm P}}

\DeclareMathOperator{\Ks}{KS}

\newcommand{\RP}{{\mathbb R\mkern-0.5mu\mathrm P}}

\newcommand{\ant}{\mathbb{A}}

\newcommand{\simtimes}{\mathbin{\widetilde{\smash{\times}}}}

\DeclareMathOperator{\Pin}{Pin}

\DeclareMathOperator{\Sp}{Spin}

\DeclareMathOperator{\Sw}{SW}

\newcommand{\cs}{\mathbin{\#}}

\begin{document}

\author{Valentina Bais and Rafael Torres}

\title[A mechanism to construct new exotic four-manifolds.]{A cut-and-paste mechanism to introduce fundamental group and construct new four-manifolds.}

\address{Scuola Internazionale Superiori di Studi Avanzati (SISSA)\\ Via Bonomea 265\\34136\\Trieste\\Italy}

\email{\{vbais, rtorres\}@sissa.it}


\maketitle

\emph{Abstract}: We introduce a simple cut-and-paste mechanism to construct both orientable and nonorientable four-manifolds from a given initial one. This mechanism alters the fundamental group while preserving other essential topological invariants. It avoids codimension two cut-and-paste fundamental group computations and fast tracks the search for fixed-point free involutions. The mechanism proves useful to unveil novel exotic irreducible smooth structures on closed four-manifolds with finite cyclic fundamental group, which include $\Q$-homology real projective four-spaces, as well as non-smoothable examples.

\section{Introduction.}\label{Introduction}

The study of existence of nondiffeomorphic smooth structures is a central theme in four-manifold theory, where the development of novel production techniques that pair well with a given diffeomorphism invariant to manufacture new examples is of pre-eminent importance. The main contribution of this paper is to record a simple construction mechanism and to sample its efficiency by building new irreducible smooth structures on closed four-manifolds. The procedure consists of four steps. 
\begin{enumerate}[label=\Roman*]
\item. Start with a pair $(M_1, M_2)$ of smooth four-manifolds that are homeomorphic, but not diffeomorphic. Remove the tubular neighborhood $\nu(\alpha_i)$ of a simple loop $\alpha_i\subset M_i$ to produce a compact four-manifold
\begin{equation*}
M_i^0\mathrel{\mathop:}= M_i\setminus \nu(\alpha_i)
\end{equation*}
for $i = 1, 2$.

\item. Build a compact oriented four-manifold $S_p$ that is a $\Q$-homology $S^1\times D^3$ with boundary $\partial S_p = S^1\times S^2$ and a compact nonorientable four-manifold $N_p$ that is a $\Q$-homology $D^3\simtimes S^1$ (the nonorientable three-disk bundle over the circle) with boundary $\partial N_p = S^2\simtimes S^1$ (the nonorientable two-sphere bundle over the circle).  A detailed description of the compact four-manifolds $S_p$ and $N_p$ is given in Section \ref{Section Oriented Trick} and Section \ref{Section nonorientable trick}, respectively. Both of these building blocks have the property that any self-diffeomorphism of their boundary extends to a self-diffeomorphism of the whole four-manifold; see Theorem \ref{Theorem LauPoe} and Theorem \ref{Theorem CMN}.

\begin{thm}\label{Theorem Extension}Let $p\in \Z_{> 0}$ be odd. 

$\bullet$ For any self-diffeomorphism $h: S^1\times S^2\rightarrow S^1\times S^2$, there is a diffeomorphism $H_p: S_p\rightarrow S_p$ such that $H_p|_{\partial} = h$. 

$\bullet$ For any self-diffeomorphism $f: S^2\simtimes S^1\rightarrow S^2\simtimes S^1$, there is a diffeomorphism $F_p: N_p\rightarrow N_p$ such that $F_p|_{\partial} = f$. 
\end{thm}

Theorem \ref{Theorem Extension} is a generalization of well-known foundational results of Laudenbach-Po\'enaru \cite{[LaudenbachPoenaru]}, C\'esar de S\'a \cite{[CesardeSa]} and Miller-Naylor \cite{[MillerNaylor]}, and it might be of independent interest.

\item. For $i=1,2$, assemble new four-manifolds
\[M_i(p)\mathrel{\mathop:}= M_i^0\cup S_p\]
if $\alpha_i\subset M_i$ is an orientation-preserving loop or
\[M_i(p)\mathrel{\mathop:}= M_i^0\cup N_p\]
if it is an orientation-reversing one. This cut-and-paste procedure introduces fundamental group, in the sense that \[\pi_1(M_i(p))\neq \pi_1(M_i),\] and its codimension three nature avoids cumbersome fundamental group computations that are unavoidable in the codimension two case: see \cite{[BaldridgeKirk1]}. The property of the building blocks $S_p$ and $N_p$ mentioned in Theorem \ref{Theorem Extension} has the following useful consequence.

\begin{thm}\label{Theorem Homeomorphism}Let $M_1$ and $M_2$ be topological four-manifolds such that there are locally flat simple loops \[\{\alpha_{M_i}\subset M_i: i = 1, 2\}\] for which there is a homeomorphism
\begin{equation*}
M_1\setminus \nu(\alpha_{M_1})\rightarrow M_2\setminus \nu(\alpha_{M_2}).
\end{equation*}
If $\alpha_{M_i}\subset M_i$ is an orientation-preserving loop, set $C_p\mathrel{\mathop:}= S_p$. Otherwise, let $C_p \mathrel{\mathop:}= N_p$. For $i=1,2$, define
\[M_i(p) \mathrel{\mathop:}=(M_i\setminus \nu(\alpha_{M_i}))\cup C_p. \]
The four-manifolds $M_1(p)$ and $M_2(p)$ are homeomorphic.
\end{thm}

\item. Find an invariant that distinguishes the diffeomorphism types of the new four-manifolds $M_1(p)$ and $M_2(p)$.
\end{enumerate}



The crux to the efficiency of the procedure is its fourth step. In principle, one would hope that the diffeomorphism invariant that discerns $M_1$ from $M_2$ also tells apart the new four-manifolds $M_1(p)$ and $M_2(p)$. This is indeed the case in all the examples under consideration in this paper, where we reckoned the Heegaard Floer homology based invariants used by Levine-Lidman-Piccirillo in \cite{[LLP]}, the Ozsv\'ath-Szab\'o invariants \cite{[OS]}, the Seiberg-Witten invariants \cite{[GompfStipsicz]}, \cite[Section 2.4]{[StipsiczSzabo2]} (as employed by Stipsicz-Szab\'o \cite{[StipsiczSzabo1], [StipsiczSzabo2]} and Baykur-Stipsicz-Szab\'o \cite{[BaykurStipsiczSzabo]}), and the $\eta$-invariant as studied by Stolz \cite{[Stolz]}. 

As a warm up result, we paint a picture of the fourth step and the utility of the mechanism through an extension of a recent result of Levine-Lidman-Piccirillo. In \cite{[LLP]}, the authors constructed the first example of an inequivalent smooth structure on a closed smooth oriented four-manifold with definite intersection form: they built a four-manifold $\mathcal{R}$ that is homeomorphic but not diffeomorphic to the connected sum $\Sigma_2\cs 4\overline{\CP^2}$, where $\Sigma_2$ is a $\Q$-homology four-sphere with $\pi_1(\Sigma_2) = \Z/2$. Using the pair $(\Sigma_2\cs 4\overline{\CP^2}, \mathcal{R})$ in the first step of our mechanism, we obtain the following generalization of their result. 

 \begin{thm}\label{Theorem A} For every odd $p \in \Z_{> 0}$, there is a closed smooth oriented four-manifold $\mathcal{R}_{2p}$ that is homeomorphic but not diffeomorphic to the connected sum $\Sigma_{2p}\cs 4\overline{\CP^2}$, where $\Sigma_{2p}$ is a $\Q$-homology four-sphere with cyclic fundamental group of order $2p$. The diffeomorphism types are distinsguished by the Ozsv\'ath-Szab\'o invariant. 
\end{thm}

To distinguish the smooth structures in \cite[Theorem 1.2]{[LLP]}, Levine-Lidman-Piccirillo show that the universal covers $\CP^2 \cs 9 \overline{\CP^2}$ and $\mathcal{E}$ of the four-manifolds $\Sigma_2 \cs 4 \overline{\CP^2}$ and $\mathcal{R}$ are non-diffeomorphic by means of their Ozsv\'ath-Szab\'o invariant \cite{[OS]}. For every odd $p\in \Z_{> 0}$, we show that the degree-two covers of our four-manifolds are obtained from $\CP^2 \cs 9 \overline{\CP^2}$ and $\mathcal{E}$, respectively, by taking a connected sum with a $\Q$-homology four-sphere. This observation, together with the properties of the Heegaard Floer smooth four-dimensional invariant, allows us to canonically adapt the proof of \cite[Theorem 1.2]{[LLP]} to distinguish such covers and, hence, also the four-manifolds $\Sigma_{2p} \cs 4 \overline{\CP^2}$ and $\mathcal{R}_{2p}$.
 


Our next theorem was inspired by recent work of Stipsicz-Szab\'o \cite{[StipsiczSzabo1], [StipsiczSzabo2]} and Baykur-Stipsicz-Szab\'o \cite{[BaykurStipsiczSzabo]}, where they built infinitely many nondiffeomorphic irreducible smooth structures on $\Sigma_{4m + 2}\cs \overline{\CP^2}$ for $\Sigma_{4m + 2}$ a $\Q$-homology four-sphere with fundamental group $\Z/(4m + 2)$. Recall that a smooth four-manifold $M$ is said to be irreducible if for every smooth connected sum decomposition $M = M_1\cs M_2$ either $M_1$ or $M_2$ is a homotopy four-sphere \cite[Definition 10.1.17]{[GompfStipsicz]}. A smooth four-manifold is minimal if it contains no smoothly embedded two-sphere of self-intersection minus one. The presence of such two-sphere in a smooth four-manifold $M$ guarantees a smooth connected sum decomposition $M = M_0\cs \overline{\CP^2}$ for some four-manifold $M_0$ \cite[Proposition 2.2.11]{[GompfStipsicz]}.

\begin{thm}\label{Theorem Main}Let $(n, p)\in \Z_{> 0}\times \Z_{\geq 3}$ and assume $p\geq 3$ is odd. Let $\Sigma_{q, k}$ be a $\Q$-homology four-sphere with finite cyclic fundamental group of order $q$ and Kirby-Siebenmann invariant $\Ks(\Sigma_{q, k}) = k$. Let $\widetilde{X}_n$ be a closed smooth oriented simply connected four-manifold that satisfies the following properties.\begin{itemize}\item There is an orientation-preserving smooth fixed-point-free involution\begin{equation}\label{Involution Hypothesis}\tau_n: \widetilde{X}_n\rightarrow \widetilde{X}_n\end{equation}and the quotient four-manifold $\widetilde{X}_n/\tau_n$ is homeomorphic to a connected sum $Y\cs \Sigma_{2, k}$, where $Y$ is a closed simply connected topological four-manifold with Kirby-Siebenmann invariant given by $\Ks(Y) = k = \Ks(\Sigma_{2p, k})$.
\item The four-manifold $\widetilde{X}_n$ has exactly two Seiberg-Witten basic classes $\pm K_n\in \mathcal{B}_{\widetilde{X}_n}$ and $SW_{\widetilde{X}_n}(\pm \mathfrak{s}_n) = \pm n^2$ for $c_1(\pm \mathfrak{s}_n) = \pm K$.\end{itemize}

Then, there exists a closed smooth oriented irreducible four-manifold $X_n(2p)$ with finite fundamental group of order $2p$ that satisfies the following properties.\begin{itemize}\item There is a homeomorphism\begin{equation}\label{Homeomorphism Conclusion}X_n(2p)\rightarrow Y\cs \Sigma_{2p, k}.\end{equation}
\item The double cover $\widehat{X}_n(p)\rightarrow X_n(2p)$ corresponding to the index-two subgroup $\Z/p\subset \pi_1(X_n(2p)) = \Z/2p$ decomposes as\begin{equation}\label{Double Cover Conclusion}\widehat{X}_n(p) = \widetilde{X}_n\cs \Sigma_p.\end{equation}
\item  There are infinitely many pairwise nondiffeomorphic irreducible four-manifolds in the homeomorphism class of $Y\cs \Sigma_{2p, k}$ and these diffeomorphisms types are distinguished by the Seiberg-Witten invariant of their double covers.
\end{itemize}

If the four-manifold $\widetilde{X}_n$ does admit an involution (\ref{Involution Hypothesis}), but has more than two Seiberg-Witten basic classes, assume it is minimal and symplectic with $b_2^+(\widetilde{X}) > 1$. Then, there exists a closed smooth oriented irreducible four-manifold $X_n(2p)$ with finite fundamental group of order $2p$ for which there is a homeomorphism (\ref{Homeomorphism Conclusion}) whose double cover is (\ref{Double Cover Conclusion}).
\end{thm}

In \cite{[StipsiczSzabo1], [StipsiczSzabo2], [BaykurStipsiczSzabo]}, the Seiberg-Witten invariants are used to distinguish the diffeomorphism types of their examples in the following way. A diffeomorphism $h: M_1\rightarrow M_2$ between two closed oriented four-manifolds lifts to a diffeomorphism of the degree-two covers $\hat{h}: \widehat{M}_1\rightarrow \widehat{M}_2$, whenever the index-two subgroups of $\pi_1(M_i)$ identifying $\widehat{M}_i \rightarrow M_i$ correspond under the map induced by $h$; see Lemma \ref{Lemma Lifts}. Stipsicz-Szab\'o obstruct the existence of $h$ through the diffeomorphism invariance of the Seiberg-Witten function $\Sw_{\widehat{M}_i}: \Sp^{\C}(\widehat{M}_i)\rightarrow \Z$ for $i = 1, 2$. Moreover, they conclude on the irreducibility of $\widehat{M}_i$ and $M_i$ through their Seiberg-Witten basic classes. We follow the same approach to prove Theorem \ref{Theorem Main}. In our case, however, the degree-two covers do decompose smoothly as a connected sum $M\cs \Sigma_p$, and we build on arguments of Baykur-Stipsicz-Szab\'o \cite{[BaykurStipsiczSzabo]} and Park \cite{[Park]}; see Section \ref{Section SW}.

As a sample of the efficiency of our simple procedure, we construct new examples of irreducible smooth structures in the following theorem.


\begin{thm}\label{Theorem Geo}Let $(m, n, p)\in \Z_{> 0}\times \Z_{> 0}\times \Z_{> 0}$ be positive integers that satisfy $4 + 5m\geq n$, $4 + 5n\geq m$, $m = 0 \mod 2$ and $p\neq 0 \mod 2$. If $m = n$,  suppose $n > 7$. There is a closed smooth irreducible four-manifold $X_{m, n}$ with finite cyclic fundamental group of order $2p$ that is homeomorphic to the connected sum\begin{equation*}m\CP^2\#n\overline{\CP^2}\#\Sigma_{2p}.\end{equation*}
\end{thm}

The integers in the statement of Theorem \ref{Theorem Geo} decompose the second Betti number of $X_{m, n}$ into\begin{center}$b_2^+(X_{m, n}(2p)) = m$ and $b_2^-(X_{m, n}(2p)) = n$.\end{center} The case of $m\neq 0 \mod 2$ and $m \leq n$ of Theorem \ref{Theorem Geo} is well-known to hold for every finite cyclic fundamental group and for infinitely many pairwise nondiffeomorphic irreducible smooth structures by modifying the constructions of closed simply connected irreducible four-manifolds of Akhmedov-Baldridge-Baykur-Kirk-Park \cite{[ABBKP]}, Akhmedov-Park \cite{[AkhmedovPark]}, Baldridge-Kirk \cite{[BaldridgeKirk1], [BaldridgeKirk2]}, Baykur-Hamada \cite{[BaykurHamada]}, Fintushel-Stern \cite{[FintushelStern2]} and Gompf \cite{[Gompf1]} as already done in \cite{[AkhmedovPark], [BaldridgeKirk2], [BaykurStipsiczSzabo], [Torres1], [Torres2]}. The case $m = 0 \mod 2$ and $p = 1$ is due to Stipsicz-Szab\'o \cite{[StipsiczSzabo1], [StipsiczSzabo2]}, Baykur-Stipsicz-Szab\'o \cite{[BaykurStipsiczSzabo]} and Arabadji-Morgan \cite{[ArabadjiPorter]}. Our proof of Theorem \ref{Theorem Geo} is based on coupling their results with Theorem \ref{Theorem Main}, and it is presented in Section \ref{Section Proof of Theorem Geo}.

Concerning the nonorientable realm, we construct new smooth structures on $\Q$-homology real projective four-spaces and non-smoothable copies of these four-manifolds through our mechanism. This extends results of Cappell-Shaneson \cite{[CappellShaneson]}, Fintushel-Stern \cite{[FintushelStern]} and Ruberman \cite{[Ruberman]}. In order to build the novel four-manifolds, we start with either Cappell-Shaneson's exotic $\RP^4$ or the non-smoothable four-manifold $\ast \RP^4$ that is homotopy equivalent to the real projective four-space; see \cite[Chapter 14]{[Wall]}, \cite[Section 2]{[Ruberman]}. We address the fourth step of our construction mechanism by computing the $\eta$-invariant \cite{[Stolz]} of the smooth four-manifolds involved in Section \ref{Section Computations Eta}. Our main result in this case is the following theorem.

\begin{thm}\label{Theorem D}For every odd $p\in \Z_{> 0}$, there are closed smooth irreducible nonorientable four-manifolds $A(2p)$ and $B(2p)$ that satisfy the following properties.\begin{enumerate}\item The Euler characteristics of these four-manifolds are equal to \[\chi(A(2p)) = 1 = \chi(B(2p))\] and their fundamental groups are isomorphic to \[\pi_1(A(2p)) = \Z/2p = \pi_1(B(2p)).\]
\item There is a homeomorphism $A(2p)\rightarrow B(2p)$.
\item The orientation 2-cover $\widehat{A}(2p)$ is diffeomorphic to the orientation 2-cover $\widehat{B}(2p)$.
\item The universal covers $\widetilde{A}(2p)$ and $\widetilde{B}(2p)$ are diffeomorphic to the connected sum $(2p - 1)(S^2\times S^2)$. 
\item There is no diffeomorphism\begin{equation*}A(2p)\cs (k - 1)(S^2\times S^2)\rightarrow B(2p)\cs (k - 1)(S^2\times S^2)\end{equation*}for any $k\in \Z_{> 0}$. 
\item The connected sums $A(2p)\cs \CP^2$ and $B(2p)\cs \CP^2$ are diffeomorphic.
\end{enumerate}

Moreover, there is a non-smoothable closed topological four-manifold $\ast A(2p)$ that is homotopy equivalent to $A(2p)$ and whose universal cover is homeomorphic to $(p - 1)(S^2\times S^2)$. 
\end{thm}

Theorem \ref{Theorem D} provides a pair of smooth orientation-reversing free $\Z/2p$-actions on a connected sum $(p - 1)(S^2\times S^2)$ that are equivariantly homeomorphic, but not equivariantly diffeomorphic. This extends work of Fintushel-Stern \cite{[FintushelStern]}, and Cappell-Shaneson \cite{[CappellShaneson]} and Gompf \cite{[Gompf]}. The last clause of Theorem \ref{Theorem D} gives examples of topological orientation-reversing free $\Z/2p$-actions on $(p - 1)(S^2\times S^2)$ whose orbit space is non-smoothable; see Hambleton-Kreck-Teichner \cite{[HambletonKreckTeichner]}. Moreover, we have collected quintessential properties of the known inequivalent smooth structures on nonorientable four-manifolds in the statement of Theorem \ref{Theorem D} for the convenience of the reader. 

The paper is organized as follows. Section \ref{Section Orientable} contains a detailed description of the building block $S_p$ and the construction mechanism in the orientable case as well as their basic properties. This includes the proofs of Theorem \ref{Theorem Extension} and Theorem \ref{Theorem Homeomorphism} in the oriented case.  A proof of Theorem \ref{Theorem A} is found in Section \ref{Section Examples Theorem A}, while Section \ref{Section Theorem Main} accommodates a proof of Theorem \ref{Theorem Main}. The proof of Theorem \ref{Theorem Geo} is discussed in Section \ref{Section Proof of Theorem Geo}. Section \ref{Section Nonorientable} is devoted to the study of the construction mechanism for nonorientable four-manifolds. In Section \ref{Section nonorientable trick}, we give a detailed description of the building block $N_p$ and the construction mechanism in the nonorientable case along with their basic properties. This includes the proofs of Theorem \ref{Theorem Extension} and Theorem \ref{Theorem Homeomorphism} in the nonorientable case. Theorem \ref{Theorem D} is proven in Section \ref{Section Proof Theorem E}.

We follow the following conventions and notation in the sequel. All our maps and manifolds are smooth unless it is stated otherwise. A pair $(X, \mathfrak{t})$ given by a four-manifold $X$ with a fixed $\Sp^\C$-structure $\mathfrak{t}$ will be called a $\Sp^\C$-structure.

\subsection{Acknowledgements} We thank Stefan Bauer, Bob Gompf, Paul Kirk, Tye Lidman, Porter Morgan, Lisa Piccirillo and Zoltan Szab\'o for useful e-mail exchanges/conversations. 


\section{Constructions of oriented four-manifolds.}\label{Section Orientable}

\subsection{Cut-and-paste procedure: oriented setting}\label{Section Oriented Trick}We construct the compact oriented four-manifold $S_p$ of the second step of our procedure as follows. Let $\alpha\subset S^1\times S^3$ be a simple loop whose homotopy class generates the infinite cyclic group $\pi_1(S^1\times S^3) = \Z$. Denote by $\alpha_p\subset S^1\times S^3$ the simple loop that represents the homology class $[\alpha]^p\in H_1(S^1\times S^3; \Z) = \Z$. The tubular neighborhood of this codimension three submanifold is $\nu(\alpha_p) = S^1\times D^3$. The main building block for our cut-and-paste construction in the orientable setting is the compact four-manifold\begin{equation}\label{bb}
S_p=(S^1 \times S^3) \setminus \nu(\alpha_p)
\end{equation}with boundary $\partial S_p = S^1\times S^2$.

\begin{example}\label{Example QSphere}The closed oriented four-manifold\begin{equation}\label{Decomposition QSphere}\Sigma_p = (D^2\times S^2)\cup S_p\end{equation}is a $\Q$-homology four-sphere whose fundamental group is finite cyclic of order $p\in \Z_{> 0}$, i.e. $\pi_1(\Sigma_p) = \Z/p$.
\end{example}

A description of the cyclic covers of the building block $S_p$ is given in our next lemma, which is of great utility in the sequel.

\begin{lemma}\label{lcover}Let $S_p$ be the four-manifold with boundary defined in (\ref{bb}). 

$\bullet$ For every $p \in \Z_{>0}$, there is a double cover
\[D_p \xrightarrow{2:1} S_p\] that corresponds to the index-two subgroup $2\Z \subset \Z$. 

$\bullet$ If $p$ is odd, then there is a diffeomorphism  $D_p \rightarrow S_p$.

$\bullet$ If $p=2q$, then there is a diffeomorphism\[D_p \rightarrow (S^1 \times S^3) \setminus \nu(\alpha_{q} \sqcup \alpha_{q}')\]
where $\alpha_{q}'$ is a parallel copy of $\alpha_{q}$.
\end{lemma}
\begin{proof}
    For any $p \in \Z_{> 0}$, the universal cover of $S_p$ is
    \[\widetilde{S}_p=\R \times (S^3 \setminus \{B_1, \dots, B_p\})\]
    where $\{B_1, \dots, B_p\} \subset S^3$ are disjointly embedded smooth three-balls. The $\Z$-action on $\widetilde{S}_p$ is generated by the map\begin{equation}\label{Phi Lemma}\phi: \widetilde{S}_p \rightarrow \widetilde{S}_p\end{equation}given by\begin{equation}(t,x) \mapsto (t+1,f_p(x)),\end{equation}
where $f_p:(S^3 \setminus \overset{p}{\underset{i = 1}{\sqcup}} B_i)\rightarrow (S^3 \setminus \overset{p}{\underset{i = 1}{\sqcup}} B_i)$ is the restriction of a self-diffeomorphism of the three-sphere that cyclically permutes the three-balls $\{B_1, \dots, B_p\}$. The degree-two cover of $S_p$ associated to the subgroup $2\Z \subset\Z= \pi_1(S_p)$ is given by the quotient of $\widetilde{S}_p$ by the $\Z$-action generated by the map $\phi^2$. The conclusion follows from the fact that the square of a cyclic permutation of odd length is still cyclic, while the square of a cyclic permutation of even length factors as a product of two cycles of half the original length.
\end{proof}

\begin{remark}\label{positionb_i's}In the following, it will be useful to consider the genus one Heegaard splitting of the three-sphere
\begin{equation}\label{heegaardS3}
S^3= (S^1 \times D^2) \cup (D^2 \times S^1).
\end{equation}We can assume that the map $f_p$ in the proof of Lemma \ref{lcover} is a rotation of $S^3$ of angle $2 \pi/p$, which restricts to a rotation along the $D^2$-factor on the solid torus on the left side of (\ref{heegaardS3}) and to a rotation along the $S^1$-factor on the solid torus on the right side. Moreover, we can choose the three-balls $\{B_1, \dots, B_p\}$ to be symmetrically embedded in the solid torus on the right side of (\ref{heegaardS3}), where they are permuted by the action of $f_p$.\end{remark}

We now prove Theorem \ref{Theorem Extension} in the orientable case. As it has been already mentioned, this theorem is an extension of a foundational result of Laudenbach-Poenar\'u \cite{[LaudenbachPoenaru]}, which states that any self-diffeomorphism of $S^1\times S^2$ extends to a self-diffeomorphism of $S^1\times D^3$. 

\begin{theorem}\label{Theorem LauPoe}(cf. Theorem \ref{Theorem Extension}). For every odd integer $p\geq 1$ and for any self-diffeomorphism\begin{equation*}h: S^1\times S^2\rightarrow S^1\times S^2,\end{equation*} there is a diffeomorphism\begin{equation*}H_p: S_p\rightarrow S_p\end{equation*} such that $H_p|_{\partial} = h$. 
\end{theorem}

We denote by $\ant: S^2\rightarrow S^2$ the antipodal involution.

\begin{proof}A result of Gluck \cite{[Gluck]} states that the diffeotopy group of $S^1 \times S^2$ is isomorphic to $\Z/2\times \Z/2\times \Z/2$. This group is generated by the reflections $(\theta, x)\mapsto (\theta, \ant(x))$ and $(\theta, x)\mapsto (-\theta, x)$ along with the self-diffeomorphism
\begin{equation}\label{gluck}
    \varphi: S^1 \times S^2 \rightarrow S^1 \times S^2
\end{equation}given by $(\theta, x) \mapsto (\theta, r_{\theta}(x))$, where $r_{\theta}: S^2\rightarrow S^2$ is a rotation of an angle $\theta$ around a fixed axis \cite{[Gluck]}, \cite[p. 156]{[GompfStipsicz]}. Such map is known as a Gluck twist. It is straight-forward to see that Theorem \ref{Theorem LauPoe} holds for the first two generators of the diffeotopy group of $S^1\times S^2$. We now show that the diffeomorphism (\ref{gluck}) admits an extension $\Phi_p: S_p \rightarrow S_p$ for any odd $p \in \Z_{> 0}$. We explicitly build the diffeomorphism $\Phi_p$ by using the identification\begin{equation}\label{Identification Sp}S_p \cong (\R \times (S^3 \setminus \overset{p}{\underset{i = 1}{\sqcup}} B_i)) / \sim_{\phi}\end{equation}coming from Lemma \ref{lcover} as follows. Without loss of generality, we work in the setting of Remark \ref{positionb_i's}: the deck transformation $\phi$ is defined in (\ref{Phi Lemma}) and the collection of three-balls $\{B_1, \dots, B_p\} \subset D^2 \times S^1$ are defined as in the proof of Lemma \ref{lcover}. 

Consider the map \[F: \R \times (S^3 \setminus \overset{p}{\underset{i = 1}{\sqcup}} B_i)  \rightarrow S^3 \setminus \overset{p}{\underset{i = 1}{\sqcup}} B_i\]defined by
\[(t,x) \mapsto r_t(x),\]where $r_t: S^3\rightarrow S^3$ is a rotation of an angle $2 \pi t$ that restricts to a rotation along the $D^2$-factor on the solid torus $D^2 \times S^1$ of the Heegaard decomposition of $S^3$ as described in Remark \ref{positionb_i's}. This solid torus contains the three-balls $\{B_1, \dots, B_p\}$. Notice that for every $t \in [0,1]$, the map $F(t, \cdot)$ is an $f_p$-equivariant self-diffeomorphism of $S^3\setminus \overset{p}{\underset{i = 1}{\sqcup}} B_i$, where $f_p$ can be chosen to act on the genus one Heegaard splitting of $S^3$ described in Remark \ref{positionb_i's}. We can now define a self-diffeomorphism\begin{equation}\label{Map Phi}\Phi: \R \times (S^3 \setminus \overset{p}{\underset{i = 1}{\sqcup}} B_i) \rightarrow \R \times (S^3 \setminus \overset{p}{\underset{i = 1}{\sqcup}} B_i)\end{equation}by\begin{equation}(t,x) \mapsto (t,F(t,x)).\end{equation}It is immediate to check that the diffeomorphism (\ref{Map Phi}) is $\phi$-equivariant. Therefore, it induces a self-diffeomorphism $\Phi_p: S_p\rightarrow S_p$. In particular, the restriction of $\Phi_p$ to $\partial S_p=S^1 \times S^2$ is a composition of $p$ Gluck twists, which is again a Gluck twist since $p$ is odd. This concludes the proof.
\end{proof}

Our construction and Theorem \ref{Theorem LauPoe} also hold in the topological category. The reader is directed to \cite[Chapter 5]{[FNOP]} regarding existence and uniqueness results of the tubular neighborhood of a locally flat simple loop in a topological four-manifold. As mentioned in the first step of our introduction, remove the tubular neighborhood of a locally flat framed simple loop $\alpha_M\subset M$ from a closed oriented topological four-manifold $M$ to obtain a compact oriented topological four-manifold $M\setminus \nu(\alpha_M)$ with boundary $\partial(M\setminus \nu(\alpha_M)) = S^1\times S^2$. Cap off its boundary with the compact topological four-manifold defined in (\ref{bb}) to construct a closed oriented topological four-manifold\begin{equation}\label{New Manifold Oriented}M(p)\mathrel{\mathop:}= (M\setminus \nu(\alpha_M)) \cup S_p.\end{equation}Theorem \ref{Theorem LauPoe} has the orientable case of Theorem \ref{Theorem Homeomorphism} as a consequence, which guarantees that the third step described in our introduction goes through.

\begin{corollary}\label{Corollary Homeomorphism Orientable}(cf. Theorem \ref{Theorem Homeomorphism}). Let $M_1$ and $M_2$ be topological four-manifolds and let $\alpha_{M_i}\subset M_i$ be a locally flat orientation-preserving simple framed loop for $i = 1, 2$. Suppose there is a homeomorphism\begin{equation*}f:M_1\setminus \nu(\alpha_{M_1})\rightarrow M_2\setminus \nu(\alpha_{M_2}).\end{equation*}The topological four-manifolds
\[M_1(p): = (M_1\setminus \nu(\alpha_{M_1}))\cup S_p \quad \text{and} \quad M_2(p)\mathrel{\mathop:}= (M_2\setminus \nu(\alpha_{M_2}))\cup S_p\]
are homeomorphic.

\end{corollary}

Independently of Corollary \ref{Corollary Homeomorphism Orientable}, the simplicity of the construction makes it straight-forward to compute several invariants of the new four-manifold (\ref{New Manifold Oriented}) as we now sample.

\begin{lemma}\label{Lemma pi orientable}The closed oriented four-manifold $M(p)$ defined in (\ref{New Manifold Oriented}) has the following basic topological invariants for every $p\in \Z_{> 0}$.
\begin{enumerate}\item The Euler characteristic is $\chi(M(p)) = \chi(M)$.
\item The signature is $\sigma(M(p)) = \sigma(M)$.
\item The equivariant intersection forms $\widetilde{I}_M$ and $\widetilde{I}_{M(p)}$ are isomorphic. 
\item The second Stiefel-Whitney class is $w_2(M(p)) = w_2(M)$.
\item If the homotopy class of the simple loop $\alpha_M\subset M$ generates the cyclic fundamental group $\pi_1(M) = \Z/q$, then the fundamental group of $M(p)$ is $\pi_1(M(p)) = \langle t: t^{q\cdot p} = 1\rangle$.
\end{enumerate}
\end{lemma}

\begin{proof}We prove the third item of the lemma. Let $\alpha_M \subset M$ be a simple loop as in (\ref{New Manifold Oriented}) and denote by $\widetilde{\alpha}_M \subset \widetilde{M}$ its preimage in the universal cover $\widetilde{M}$ of $M$. The universal cover $\widetilde{M}(p)$ of $M(p)$ is constructed by gluing together copies of $\widetilde{M} \setminus \nu(\widetilde{\alpha}_M)$ and $\widetilde{S}_p/ \sim_{\phi^k}$ by matching each boundary component of $\widetilde{M} \setminus \nu(\widetilde{\alpha}_M)$ with a boundary component of $\widetilde{S}_p/ \sim_{\phi^k}$, where $\widetilde{S}_p$ is the universal cover of $S_p$ (see Lemma \ref{lcover}) and $k$ is the order of the homotopy class $[\alpha_M] \in \pi_1(M)$. In particular, if $k = \infty$ we set $\widetilde{S}_p/\sim_{\phi^{\infty}}=\widetilde{S}_p $. It follows that the equivariant intersection form $\widetilde{I}_{M(p)}$ of $M(p)$ is
\[\widetilde{I}_{M(p)} = \widetilde{I}_M \otimes_{\Z[\pi_1(M)]} \Z[\pi_1(M(p))].\] A proof of the other items of Lemma \ref{Lemma pi orientable} is obtained from straight-forward homological computations, Novikov's additivity \cite[Remark 9.1.7]{[GompfStipsicz]} and the Seifert-van Kampen theorem. 
\end{proof}

\begin{lemma}\label{Lemma Decomposition}Let $M$ be a closed oriented four-manifold and let $\alpha_M\subset M$ be a simple framed loop whose homotopy class is $[\alpha_M] = 0\in \pi_1(M)$. There is a diffeomorphism\begin{equation}\label{Diffeomorphism CS Null}M(p) = (M\setminus \nu(\alpha_M))\cup S_p \rightarrow M \cs \Sigma_p ,\end{equation} where $\Sigma_p$ is a  $\Q$-homology four-sphere with finite cyclic fundamental group of order $p$.
\end{lemma}

\begin{proof}The existence of the diffeomorphism (\ref{Diffeomorphism CS Null}) can be seen from the decompositions (\ref{Decomposition QSphere}) in Example \ref{Example QSphere} and\begin{equation}M = M\cs S^4 = M\cs ((D^2\times S^2)\cup (S^1\times D^3)).\end{equation}\end{proof}

Lemma \ref{lcover} and Lemma \ref{Lemma Decomposition} have the following consequence.

\begin{lemma}\label{Lemma Decomposition Two Cover}Let $M(2)$ be a closed oriented four-manifold whose universal cover is $\widetilde{M}$ and let $\alpha\subset M(2)$ be a closed simple framed loop such that its homotopy class generates the group $\pi_1(M(2)) = \Z/2$. If $p\in \Z_{> 0}$ is odd, then the degree-two cover $\widehat{M}(p)$ of the closed oriented four-manifold\begin{equation*}M(2p) = (M(2)\setminus \nu(\alpha))\cup S_p\end{equation*} is diffeomorphic to\begin{equation*}\widehat{M}(p) = \widetilde{M}\cs \Sigma_p.\end{equation*} 
\end{lemma}

We finish this section with the following observation regarding lifts of diffeomorphisms to covers, whose proof is left to the reader.

\begin{lemma}\label{Lemma Lifts} Let\begin{equation}\label{Diffeo Down}f:M_1\rightarrow M_2\end{equation}be a diffeomorphism between two four-manifolds $M_1$ and $M_2$ with fundamental group $\pi_1(M_i) = \Z/2p$ and let $\widetilde{M}_i$ be the degree-two cover of $M_i$ for $i = 1, 2$. The diffeomorphism (\ref{Diffeo Down}) lifts to a diffeomorphism \[\tilde{f}: \widetilde{M}_1\rightarrow \widetilde{M}_2.\]
\end{lemma}

Although Lemma \ref{Lemma Lifts} clearly holds under much less restrictive hypothesis, this phrasing suffices for the purposes in this manuscript. 

\section{Examples of Theorem \ref{Theorem A}.}\label{Section Examples Theorem A}

In order to construct the four-manifolds $\mathcal{R}_{2p}$ of Theorem \ref{Theorem A},  we make heavy use of the tools developed by Levine-Lidman-Piccirillo in \cite{[LLP]}.

\subsection{Raw materials for Theorem \ref{Theorem A}} We will employ as a raw material the exotic four-manifolds $\mathcal{E}$ and $\mathcal{R}$ constructed in \cite[Theorem 1.2]{[LLP]}.

\begin{theorem}[Levine-Lidman-Piccirillo \cite{[LLP]}]\label{Theorem LLP}\
\begin{enumerate}
    \item There exists a four-manifold $\mathcal{E}$, which is homeomorphic but not diffeomorphic to $\CP^2 \cs 9 \overline{\CP^2}$.
    \item There exists a four-manifold $\mathcal{R}$, which is homeomorphic but not diffeomorphic to $\Sigma_2 \cs 4 \overline{\CP^2}$, where $\Sigma_2$ is a $\Q$-homology four-sphere with $\pi_1(\Sigma_2)\cong \Z/2$.
\end{enumerate}
\end{theorem}

We refer the reader to \cite[Section 6.2]{[LLP]} for the precise definition of the four-manifolds $\mathcal{E}$ and $\mathcal{R}$. The four-manifolds of Theorem \ref{Theorem A} are defined by applying the construction described in Section \ref{Section Oriented Trick} to the exotic four-manifold $\mathcal{R}$ in the second item of Theorem \ref{Theorem LLP}. 

For any integer $p \geq 1$, we define the closed oriented four-manifold\begin{equation}\label{R}
    \mathcal{R}_{2p}=(\mathcal{R} \setminus \nu(\alpha_{\mathcal{R}})) \cup S_p,
\end{equation}where $\alpha_{\mathcal{R}} \subset \mathcal{R}$ is a simple closed curve whose homotopy class generates the group $\pi_1(\mathcal{R}) = \Z/2$. The homeomorphism type of each four-manifold (\ref{R}) is pinned down in the following lemma.
\begin{lemma}\label{Lemma O Homeomorphism Type}
For every integer $p \geq 1$, there is a homeomorphism\begin{equation}\label{cs}\mathcal{R}_{2p} \rightarrow\Sigma_{2p} \cs 4 \overline{\CP^2},
\end{equation}where $\Sigma_{2p}$ is a $\Q$-homology four-sphere with fundamental group $\Z/{2p}$.
\end{lemma}
\begin{proof}The $\Q$-homology four-sphere on the right side of (\ref{cs}) can be deconstructed as\begin{equation}\label{Decomposition Homology Sphere}\Sigma_{2p} = (\Sigma_2 \setminus \nu(\alpha_{\Sigma_2}))\cup S_p,\end{equation} where $\alpha_{\Sigma_2}\subset \Sigma_2$ is a simple loop whose homotopy class generates the fundamental group $\pi_1(\Sigma_2) \cong \Z/2 \cong \pi_1(\Sigma_2\cs 4 \overline{\CP^2})$. The second clause of Theorem \ref{Theorem LLP} implies that there is a homeomorphism\begin{equation*}(\Sigma_2\cs 4\overline{\CP^2})\setminus \nu(\alpha_{\Sigma_2})\rightarrow \mathcal{R} \setminus \nu(\alpha_{\mathcal{R}}).\end{equation*}Such homeomorphism can be extended to a homeomorphism as in (\ref{cs}) by the decomposition of $\Sigma_{2p}$ in (\ref{Decomposition Homology Sphere}) and Corollary \ref{Corollary Homeomorphism Orientable}. Alternatively, the existence of the homeomorphism (\ref{cs}) also follows from work of Hambleton-Kreck \cite[Theorem C]{[HambletonKreck]}.\end{proof}

To distinguish the diffeomorphism class of $\mathcal{R}_{2p}$ from the one of the connected sum $\Sigma_{2p}\cs 4\overline{\CP^2}$, we look at the double cover\[\mathcal{E}_p\xrightarrow{2:1} \mathcal{R}_{2p}\]that corresponds to the index-two subgroup $\Z/p \subset \Z/2p= \pi_1(\mathcal{R}_{2p})$. 


\begin{proposition}\label{cover}If $p \in \Z_{> 0}$ is odd, the degree-two cover of $\mathcal{R}_{2p}$ that corresponds to the index-two subgroup $\Z/p \subset \Z/2p$ is diffeomorphic to the four-manifold\[\mathcal{E}_p=(\mathcal{E} \setminus \nu(\gamma))\cup S_p = \mathcal{E}\cs \Sigma_p,\]where $\mathcal{E}$ is the closed four-manifold of the first clause of Theorem \ref{Theorem LLP}, $\gamma\subset \mathcal{E}$ is a simple loop and $S_p$ is the compact four-manifold defined in (\ref{bb}).\end{proposition}\begin{proof}The proposition follows from the definition of $\mathcal{R}$ in \cite{[LLP]} as the quotient of $\mathcal{E}$ by a fixed-point-free orientation-reversing involution, together with Lemma \ref{Lemma Decomposition Two Cover}. \end{proof}

\subsection{Diffeomorphism obstruction}\label{Section Diffeomorphism obstruction}
In light of the definition of $\mathcal{E}$ in \cite[Section 6.2]{[LLP]}, Proposition \ref{cover} implies the existence of a decomposition\begin{equation}\label{decomposition}\mathcal{E}_p = V \cup_{\zeta} V_p,\end{equation}where $V = C_0 \cup_{\phi_0 \circ \sigma} Z$ is the building block defined in \cite[Construction 6.4]{[LLP]}. The compact four-manifold $C_0$ is contractible, $Z$ is obtained by attaching five $-1$-framed $2$-handles to $\partial C_0$, the boundary $\partial V$ is the result of $0$-surgery on $S^3$ along the Stevedore knot and $\zeta$ is any of its orientation-reversing fixed-point-free involutions. On the other hand, for every $p \geq 1$ we have that\[V_p = C_0^p \cup_{\phi_0 \circ \sigma} Z\]where $C_0^p = (C_0\setminus \nu(\gamma_{C_0}))\cup S_p$ for $\gamma_{C_0}\subset C_0$ a simple loop. We assemble\[\mathcal{E}'_p=V' \cup_{\zeta} V'_p,\]where $V' = C_0 \cup_{\phi_0} Z$ and $V'_p=C_0^p \cup_{\phi_0} Z$ are made up by the same building blocks as $V$ and $V_p$ with different gluing maps. We refer the reader to \cite[Construction 6.4]{[LLP]} for further details.

The aim of this section is to provide a proof of the following proposition. 

\begin{proposition}\label{prop}
For every odd $p \in \N$, the following holds.
    \begin{enumerate}
        \item The four-manifold $\mathcal{E}'_p$ is diffeomorphic to\[X_p = ((\CP^2 \cs 9\overline{\CP^2})\setminus \nu(\gamma)) \cup S_p = (\CP^2 \cs 9\overline{\CP^2})\cs \Sigma_p,\]where $\gamma \subset \CP^2 \cs 9 \overline{\CP^2}$ is a simple loop.
        \item The four-manifold $\mathcal{E}_p$ is homeomorphic but not diffeomorphic to $X_p$.
    \end{enumerate}
\end{proposition}

We proceed to provide some technical lemmas that are needed to distinguish the smooth structures of Theorem \ref{Theorem A}. These results are generalizations of \cite[Proposition 6.9]{[LLP]} and \cite[Lemma 6.3]{[LLP]}, and make use of the Ozsváth-Szabó four-manifold invariants \cite[Section 2.4]{[OS]}. Before delving into the lemmas and for the convenience of the reader,  we briefly recall the definition of the Ozsv\'ath--Szab\'o mixed invariant for closed smooth four-manifolds. The field with two elements is denoted by $\mathbb{F}$ in the sequel.

Given a closed oriented three-manifold $Y$ with a $\Sp^{\mathbb{C}}$-structure $\mathfrak{s}$, we denote by $HF^-(Y,\mathfrak{s})$ the minus version of its Heegaard Floer homology. This is a finitely generated graded module over $\mathbb{F}[U]$, whose $U$-torsion part is denoted by $HF_{\text{red}}(Y,\mathfrak{s})$. 

\begin{example}
    The three-sphere has a unique $\Sp^{\mathbb{C}}$-structure $(S^3, \mathfrak{s})$ and \[HF^-(S^3,\mathfrak{s})\cong \mathbb{F}[U]\] where the degree of the unit element is $\text{deg}(1)=0$.
\end{example}

Given a $\Sp^{\mathbb{C}}$-cobordism $(W,\mathfrak{t})$ between $(Y_1,\mathfrak{s}_1)$ and $(Y_2,\mathfrak{s}_2)$, there is an induced map 
\[F^-_{(W,\mathfrak{t})}\colon  HF^-(Y_1,\mathfrak{s}_1) \longrightarrow HF^-(Y_2,\mathfrak{s}_2).\]Moreover, there is a canonical projection
\[\Pi\colon  HF^-(Y,\mathfrak{s}) \cong HF_{\text{red}}(Y,\mathfrak{s})\oplus \mathbb{F}[U]^k\longrightarrow HF_{\text{red}}(Y,\mathfrak{s}),\]
see \cite[Section 2.3]{[OS]}. For more details on this invariant, we refer the reader to \cite{MO}.

\begin{definition}
    Let $(W, \mathfrak{t})$ be a compact $\Sp^{\mathbb{C}}$-structure with non-empty boundary $\partial(W,\mathfrak{t})=(Y,\mathfrak{s})$, where $\mathfrak{s}=\mathfrak{t}|_Y$ is non-torsion (i.e. its first Chern class $c_1(\mathfrak{s})$ is not a torsion element in $ H_2(Y;\Z)$). The relative invariant of $(W, \mathfrak{t})$ is 
    \[\Psi_{W,\mathfrak{t}}= \Pi (F^-_{(W \setminus D^4, \mathfrak{t}|)}(1))\in HF^-_{\text{red}}(Y, \mathfrak{s}).\]   
\end{definition} 

The duality between $HF_{\text{red}}(Y, \mathfrak{s})$ and $HF_{\text{red}}(\overline Y, \mathfrak{s})$  \cite[page 376]{[OS2]} allows us to define a pairing
\[\langle \cdot, \cdot \rangle_Y \colon HF_{\text{red}}(Y, \mathfrak{s}) \otimes HF_{\text{red}}(\overline Y, \mathfrak{s}) \longrightarrow \mathbb{F}.\]With this, we are ready to define the Ozsv\'ath--Szab\'o mixed invariant for smooth four-manifolds.

Let $X$ be a closed oriented four-manifold with $b_2^+(X)=1$ and indefinite intersection form, and let $L\subseteq H_2(X;\mathbb{Q})$ be a one-dimensional subspace on which the intersection form vanishes. We say that a $\Sp^{\mathbb{C}}$-structure $(X, \mathfrak{t})$ is allowable if $c_1(\mathfrak{t}|_L)$ is non-vanishing. An embedded, oriented, separating three-manifold $Y\subseteq X$ is an admissible cut if the image in $H_2(X;\mathbb{Q})$ of the map induced by the inclusion $Y \rightarrow X$ is equal to $L$. At this point, we can write
\[X=W_1 \cup_Y W_2\]
where $W_1$ and $W_2$ are oriented in such a way that $\partial W_1=Y = \overline{\partial W_2}$. 

\begin{definition}
    The closed four-manifold invariant of the $\Sp^{\mathbb{C}}$-structure $(X, \mathfrak{s})$ with respect to the one-dimensional subspace $L\subseteq H_2(X; \mathbb{Q})$ is
    \[\Phi_{X,\mathfrak{t},L}=\langle \Psi_{(W_1, \mathfrak{t}_1)}, \Psi_{(W_2, \mathfrak{t}_2)}\rangle_Y, \]
where $\mathfrak{t}_1=\mathfrak{t}|_{W_1}$, $\mathfrak{t}_2=\mathfrak{t}|_{W_2}$ and $\mathfrak{s}=\mathfrak{t}|_Y$. 
\end{definition}

We refer the reader to \cite{[OS]} and \cite[Section 6.1]{[LLP]} for further background on these invariants. We now go into the aforementioned lemmas that are needed to obstruct the existence of the diffeomorphisms of Theorem \ref{Theorem A} and distinguish the smooth structures.


\begin{lemma}\label{l1}
    For any $p \geq 1$, there is a $\Sp^{\C}$-structure $(V_p, \mathfrak{t}_p)$ such that $\Psi_{V_p, \mathfrak{t}_p} \neq 0$. Moreover, all such ${\Sp}^{\C}$-structures can be chosen to coincide on $Z \subset V_p$. In particular, $V_p$ is homeomorphic but not diffeomorphic to $X_0(Q_0)_p' \cs 4 \overline{\CP^2}$ for\begin{equation*}X_0(Q_0)_p'=(X_0(Q_0) \setminus \nu(\gamma)) \cup S_p,\end{equation*} where $X_0(Q_0)$ is the $0$-trace of the Stevedore knot and $\gamma \subset X_0(Q_0)$ is a null-homotopic simple loop.
\end{lemma}
\begin{proof}
    The case $p=1$ is \cite[Proposition 6.9]{[LLP]}. Since $C_0^p \setminus B^4$ is a rational homology cobordism between $S^3$ and $\partial C_0^p$, the same proof holds for $p>1$ by replacing $C_0$ with $C_0^p$ in the argument. 
\end{proof}

\begin{lemma}\label{l2}
    Let $K$ be a fibered amphichiral knot of genus $g$. Let $W_1$ and $W_2$ be compact oriented 4-manifolds with positive second Betti number and boundary $\partial W_1 = S^3_0(K) = \partial W_2$. Suppose there are ${\Sp}^{\C}$ structures $\{(W_i, \mathfrak{u}_i): i = 1, 2\}$ such that $c_1(\mathfrak{u}_i)$ evaluates to $2g-2$ on a generator of $H_2(S^3_0(K))$ and for which $\Psi_{W_i,\mathfrak{u}_i} \neq 0$. Let $X$ be the closed four-manifold that is obtained by gluing together $W_1$ and $W_2$ via any choice of orientation-reversing diffeomorphism of $S^3_0(K)$, and let $L \subset H_2(X;\mathbb{Q})$ be the line corresponding to the inclusion\begin{equation*}H_2(S^3_0(K)) \rightarrow H_2(X).\end{equation*} Then, there is a ${\Sp}^{\C}$-structure $(X, \mathfrak{t})$ for which $\Phi_{X, \mathfrak{t},L}\neq 0$.
\end{lemma}

The proof of Lemma \ref{l2} is a straightforward adaptation of the one of \cite[Lemma 6.3]{[LLP]}. We now use the previous pair of lemmas to prove Proposition \ref{prop}.

\begin{proof}[Proof of Proposition \ref{prop}]
Item (1) follows from the construction and the fact that $\mathcal{E}'$ is an exotic $\CP^2 \cs 9 \overline{\CP^2}$ by \cite[Theorem 6.7]{[LLP]}. By Lemma \ref{Lemma Decomposition}, there is a degree-two cover 
\[X_p \cong\mathcal{E}'_p \xrightarrow{2:1} \Sigma_{2p} \cs 4 \overline{\CP^2}. \]
Let us now prove Item (2). In particular, $\mathcal{E}_p$ and $\mathcal{E}'_p$ are homeomorphic by construction, since $\mathcal{E}$ and $\mathcal{E}'$ are homeomorphic by \cite[Theorem 6.7]{[LLP]}. To conclude the proof, we need to show that the degree-two covers $\mathcal{E}_p\rightarrow \mathcal{R}_{2p}$ and $\mathcal{E}'_p\rightarrow\Sigma_2 \cs 4 \overline{\CP^2}$ are not diffeomorphic. Consider the decomposition
\begin{equation}
    \mathcal{E}_p=V_a \cup T_{Q_0} \cup V_b^p
\end{equation}
where $V_a$ and $V_b^p$ are copies of $V$ and $V_p$ respectively without the $T^-_{Q_0}$ component, see \cite[Construction 6.4]{[LLP]}. As for the case of $\mathcal{E}$, we see that $H_2(\mathcal{E}_p)$ has a basis
\[(E_1, \dots, E_4, E_5^p, \dots, E_8^p, \Xi, \Theta)\]
where:
\begin{itemize}
\item $(E_1, \dots, E_4)$ is a diagonal basis for $V_a$:
\item $(E_5^p, \dots, E_8^p$ is a diagonal basis for $V_b^p$;
\item $(\Xi, \Theta)$ is the basis for $H_2(T_{Q_0})$ from \cite[Remark 3.11]{[LLP]}. In particular, $\Xi$ is represented by a capped-off Seifert surface for $Q_0$ in $S^3(Q_0)$ and $\Theta$ is represented by a $2$-sphere of self-intersection $-2$ and $\Xi \cdot \Theta=1$.
\end{itemize}

Lemma \ref{l1} and Lemma \ref{l2} imply that there is a ${\Sp}^{\C}$-structure $(\mathcal{E}_p, \mathfrak{t})$ such that $\Phi_{\mathcal{E}_p,\mathfrak{t},L} \neq 0$, where $L = \text{Span}(\Xi) \subset H_2(\mathcal{E}_p;\mathbb{Q})$. By \cite[Lemma 6.2]{[LLP]}, we have that $\Xi$ cannot be represented by a surface of genus less than two. If there were a diffeomorphism\[\mathcal{E}_p \rightarrow X_p,\]then one could build a spherical representative for $\Xi$ as in the proof of \cite[Theorem 6.7]{[LLP]}. This would lead to a contradiction.\end{proof}

\begin{remark}If $p$ is even, then the four-manifold $D_p$ of Lemma \ref{lcover} has two boundary components and the double cover $\mathcal{E}_p$ of $\mathcal{R}_{2p}$ is\begin{equation*}\mathcal{E}_p=(\mathcal{R}\setminus \nu(\alpha_{\mathcal{R}})) \cup D_p \cup (\mathcal{R} \setminus \nu(\alpha_{\mathcal{R}})).\end{equation*}We have not been able to distinguish the smooth structure $\mathcal{R}_{2p}$ using the Heegaard Floer invariants of $\mathcal{E}_p$ as in \cite{[LLP]}. 
\end{remark}

\subsection{Proof of Theorem \ref{Theorem A}}\label{Section Proof Theorem A}The four-manifold $\mathcal{R}_{2p}$ was constructed in (\ref{R}) for any odd $p\in \Z_{> 0}$ and its homeomorphism type was pinned down in Lemma \ref{Lemma O Homeomorphism Type}. The claim that there does not exist a diffeomorphism $\mathcal{R}_{2p}\rightarrow \Sigma_{2p} \cs 4 \overline{\CP^2}$ follows from Proposition \ref{prop}. Indeed, if any such  diffeomorphism between $\mathcal{R}_{2p}$ and $\Sigma_{2p} \cs 4 \overline{\CP^2}$ were to exist, it would induce a diffeomorphism between the degree-two covers\begin{equation*}\mathcal{E}_p\rightarrow((\CP^2 \cs 9 \overline{\CP^2})\setminus \nu(\gamma))\cup S_p = (\CP^2 \cs 9 \overline{\CP^2})\cs \Sigma_p\end{equation*}by Lemma \ref{Lemma Decomposition Two Cover} and Lemma \ref{Lemma Lifts}. However, Proposition \ref{prop} indicates that these four-manifolds are not diffeomorphic.\hfill $\square$

\section{Mechanism in Theorem \ref{Theorem Main}.}\label{Section Theorem Main}

\subsection{Seiberg-Witten invariants of smoothly reducible double covers}\label{Section SW}Let $\Sp^\C(X)$ be the set of $\Sp^\C$-structures on a closed oriented 4-manifold $X$. We regard the Seiberg-Witten invariant of $X$ as a function\begin{equation*}\Sw_X:\Sp^\C(X)\rightarrow \Z.\end{equation*}A cohomology class $K\in H^2(X; \Z)$ is a Seiberg-Witten basic class if $\Sw_X(\mathfrak{s})\neq 0$ for a $\Sp^\C$-structure $(X, \mathfrak{s})$ with $c_1((X, \mathfrak{s}) = K$. The set of Seiberg-Witten basic classes of the four-manifold $X$ is denoted by $\mathfrak{B}_X$. We direct the reader towards \cite{[BaykurStipsiczSzabo],[GompfStipsicz], [StipsiczSzabo2],[Taubes]} and the references there for the necessary background on Seiberg-Witten theory.

The proofs of our main results require us to determine the Seiberg-Witten invariants of connected sums of four-manifolds $X = X_1\cs X_2$. These invariants are known to be determined in terms of the summands $X_1$ and $X_2$ of the connected sum $X$. We start this section by recalling a pair of results for such purpose. For any $\Sp^\C$-structure $(X_1, \mathfrak{s}_{X_1})$, denote its extension to $X = X_1\cs X_2$ by $(X, \mathfrak{s}_X)$.  Fintushel-Stern's blow up formula \cite{[FintushelStern1]} establishes a relation between the invariants of $X$ and $X_1$.

\begin{proposition}\label{Proposition BU Formula} \cite{[FintushelStern1]}; cf. \cite[Proposition 2.8]{[BaykurStipsiczSzabo]}. Let $\{(X_i, \mathfrak{s}_i): i = 1, 2\}$ be $\Sp^\C$-structures where the underlying four-manifolds are closed and oriented. Consider nonnegative integers $\{n_i: i = 1, \ldots, k\}$ given by $\langle c_1(\mathfrak{s}_2), e_i\rangle = 2n_i + 1$ that satisfy the identity\begin{equation*}\frac{1}{4}(c_1(\mathfrak{s}_1^2 - 2\chi(X_1) - 3\sigma(X_1)) - \overset{k}{\underset{i = 1}\sum} n_i(n_i + 1)\geq 0.\end{equation*}The Seiberg-Witten invariants satisfy\begin{equation*}\Sw_{X_1}(\mathfrak{s}_1) = \pm \Sw_{X_1\cs X_2}(\mathfrak{s}_1\cs \mathfrak{s}_2).\end{equation*}\end{proposition}

Baykur-Stipsicz-Szab\'o observed that the next proposition is a special case of the Seiberg-Witten blow up formula in \cite[Corollary 2.11]{[BaykurStipsiczSzabo]}. The result was also obtained by Kotschick-Morgan-Taubes \cite[Proposition 2]{[KotschickMorganTaubes]}.

\begin{proposition}\label{Proposition SW Connected Sum}(\cite{[BaykurStipsiczSzabo], [KotschickMorganTaubes]}). Let $X_1$ and $X_2$ be closed oriented four-manifolds and suppose that $b_1(X_2) = 0 = b_2^+(X_2)$. For any $\Sp^\C$-structure $(X_1, \mathfrak{s}_{X_1})$, denote by $(X, \mathfrak{s}_X)$its extension to the connected sum $X = X_1\cs X_2$. The corresponding Seiberg-Witten invariants satisfy\begin{equation*}\Sw_{X}(\mathfrak{s}_X) = \pm \Sw_{X_1}(\mathfrak{s}_{X_1}).\end{equation*}In particular, if $X_1$ has nontrivial Seiberg-Witten invariant, then so does $X$.
\end{proposition}

The next proposition is the main technical contribution in this section. Besides the pair of proposition stated before, we employ arguments of Baykur-Stipsicz-Szab\'o \cite[Proof of Proposition 6.2]{[BaykurStipsiczSzabo]} and of Park \cite[Proof of Theorem 6]{[Park]} in its proof. 

\begin{proposition}\label{Proposition SW Procedure} Let $\widetilde{M}$ be a closed oriented simply connected four-manifold that satisfies one of the following two properties. 

$\bullet$ It has exactly two Seiberg-Witten basic classes $\pm K_{\widetilde{M}}\in \mathfrak{B}_{\widetilde{M}}$.

$\bullet$ It admits a symplectic structure $(\widetilde{M}, \omega_{\widetilde{M}})$ and it is minimal with $b_2^+(\widetilde{M}) > 1$. 

Suppose there is an orientation-preserving fixed-point-free involution
\begin{equation*}
\tau:\widetilde{M}\rightarrow \widetilde{M}
\end{equation*}
and denote by $M(2) = \widetilde{M}/\tau$ the quotient closed smooth oriented four-manifold. Let $\alpha_{M(2)}\subset M(2)$ a simple loop whose homotopy class generates the fundamental group $\pi_1(M(2)) = \Z/2$. 

The four-manifold\begin{equation*}M(2p) = (M(2)\setminus \nu(\alpha_{M(2)}))\cup S_p\end{equation*}is irreducible provided $p\in \Z_{> 0}$ is  odd. Furthermore, there is a $\Sp^\C$-structure $(\widehat{M}(p), \mathfrak{s}_{\widehat{M}(p)})$ on the degree-two cover\begin{equation*}\widehat{M}(p)\xrightarrow{2:1} M(2p)\end{equation*}that corresponds to the subgroup $\Z/p\subset \pi_1(M(2p)) = \Z/2p$ for which the value of the Seiberg-Witten invariant satisfies\begin{equation}\label{SW CS}\Sw_{\widehat{M}(p)}(\mathfrak{s}_{\widehat{M}(p)}) = \Sw_{\widetilde{M}\cs \Sigma_p}(\mathfrak{s_{\cs}}) = \Sw_{\widetilde{M}}(\mathfrak{s}_{\widetilde{M}})\end{equation}for $\Sp^\C$-structures $(\widetilde{M}\cs \Sigma_p, \mathfrak{s}_{\cs})$ and $(\widetilde{M}, \mathfrak{s}_{\widetilde{M}})$.
\end{proposition}

A symplectic four-manifold is minimal if and only if it contains no embedded two-sphere of self-intersection minus one \cite[Remark 10.2.5 (a)]{[GompfStipsicz]}. The reader might have already noticed that the properties that the simply connected four-manifold $\widetilde{M}$ of Proposition \ref{Proposition SW Procedure} is required to have imply that it is irreducible \cite{[Taubes], [Kotschick], [HamiltonKotschick]}. 

\begin{proof} Proceed by contradiction and suppose that there is a smooth connected sum decomposition\begin{equation}\label{Smooth CS Decomposition}M(2p) = M_1\cs M_2,\end{equation}where neither of the four-manifolds $M_1$ and $M_2$ is a homotopy four-sphere. Since the finite cyclic group $\pi_1(M(2p)) = \Z/2p$ does not decompose as a free product of two non-trivial groups, without loss of generality we can assume that $\pi_1(M_1) = \Z/2p$, $\pi_1(M_2) = \{1\}$ and $b_2(M_2) \geq 1$. This yields a decomposition of the degree-two cover of (\ref{Smooth CS Decomposition}) as\begin{equation}\label{Decomposition Two Cover 2}\widehat{M}(p) = \widehat{M}_1\cs 2M_2,\end{equation}where $\widehat{M}_1$ is the degree-two cover of $M_1$. We also have the decomposition\begin{equation}\label{Decomposition Recall Two Cover}\widehat{M}(p)=\widetilde{M} \cs \Sigma_p\end{equation}from Lemma \ref{Lemma Decomposition Two Cover}. 

These identifications and Proposition \ref{Proposition SW Connected Sum} imply that there are bijections\begin{equation}\label{Bijections Basic Classes}\mathfrak{B}_{\widehat{M}_1\cs 2 M_2}\rightarrow \mathfrak{B}_{\widetilde{M}\cs \Sigma_p}\rightarrow \mathfrak{B}_{\widetilde{M}}\end{equation}among the Seiberg-Witten basic classes. If $\widetilde{M}$ has non-trivial Seiberg-Witten invariant, then $\widehat{M}(p)$ also has non-trivial Seiberg-Witten invariant and $b_2^+(2M_2) = 0$ \cite[Theorem 2.4.6]{[GompfStipsicz]} (this is the case in both instances of the proposition (\cite{[Taubes]})). This implies that the intersection form of $Q_{2M_2}$ is negative-definite and Donaldson's theorem \cite{[Donaldson]} says that $Q_{2M_2} \cong 2 b_2(M_2) \langle -1 \rangle$. Let $\{e^1_1, e^2_1, \ldots, e^1_{b_2(M_2)}, e^2_{b_2(M_2)}\}$ be a basis for the second cohomology group $H^2(2M_2; \Z) = \Z^{2b_2(M_2)}$ with $(e^i_j)^2 = -1$ for each $i = 1, 2$ and $j = 1, \ldots, b_2(M_2)$. 

From (\ref{Decomposition Two Cover 2}) and using the neck pinching argument \cite{[Donaldson], [Kotschick]}, we conclude that there is at least one basic class $K_{\widehat{M}_1}\in \mathcal{B}_{\widehat{M}_1}$ and the Seiberg-Witten invariant is $\Sw_{\widehat{M}_2}(\mathfrak{s}_{\widehat{M}_1})\neq 0$ for $c_1(\mathfrak{s}_{\widehat{M}_1}) = K_{\widehat{M}_1}$. The set of basic classes $\mathfrak{B}_{\widehat{M}_1\cs 2 M_2}$ consists of elements of the form\begin{equation}\label{Basic Classes Contradiction}K_{\widehat{M}_1} \pm e^1_1 \pm e^2_2 \pm \cdots \pm e^1_{b_2(M_2)} \pm e^2_{b_2(M_2)}\end{equation}for $K_{\widehat{M}_1}\in \mathfrak{B}_{\widehat{M}_1}$ \cite[Theorem 2.4.10, Remark 2.4.11]{[GompfStipsicz]}. At this point, we employ an argument due to Baykur-Stipsicz-Szab\'o \cite{[BaykurStipsiczSzabo]} to prove the first instance of Proposition \ref{Proposition SW Procedure}. Consider the map\begin{equation}\label{Map C1}c_1: \Sp^{\C}(\widehat{M}(p))\rightarrow H^2(\widehat{M}(p); \Q)\end{equation}
assigning to a given $\Sp^{\C}$-structure $(\widehat{M}(p), \mathfrak{s})$ its first Chern class with $\Q$-coefficients. The punchline is that the number of elements of the image of the set\begin{equation}\label{Set SpinC}\{\mathfrak{s}\in \Sp^\C(\widehat{M}(p)): \Sw_{\widehat{M}(p)}(\mathfrak{s})\neq 0\}\end{equation}under the map (\ref{Map C1}) differs for the decompositions (\ref{Decomposition Two Cover 2}) and (\ref{Decomposition Recall Two Cover}), and this yields a contradiction. Suppose the four-manifold $\widetilde{M}$ has exactly two Seiberg-Witten basic classes. Proposition \ref{Proposition BU Formula} and Proposition \ref{Proposition SW Connected Sum} along with the decomposition (\ref{Decomposition Recall Two Cover}) imply that the image of the set (\ref{Set SpinC}) under the map (\ref{Map C1}) has two elements, given that the four-manifold $\widetilde{M}$ is assumed to have exactly two Seiberg-Witten basic classes. On the other hand, the decomposition (\ref{Decomposition Two Cover 2}) implies that the image of the set (\ref{Set SpinC}) under (\ref{Map C1}) contains more than two elements by (\ref{Basic Classes Contradiction}). This is a contradiction and we conclude that the four-manifold $M(2p)$ is irreducible. 

We now address the second instance of the proposition using an argument due to Park \cite[Proof of Theorem 6]{[Park]} to finish the proof. Suppose that $\widetilde{M}$ admits a symplectic structure and that it is minimal. The previous arguments show that there is an isomorphism\begin{equation}\label{Isomorphism}\xi: H^2(\widetilde{M}; \Z)\rightarrow H^2(\widehat{M}_1\cs 2 M_2; \Z)\end{equation}that preserves the cup product pairing and that restricts to a bijection $\mathfrak{B}_{\widetilde{M}}\rightarrow \mathfrak{B}_{\widehat{M}_1\cs 2 M_2}$. Taubes showed that $c_1(\widetilde{M}, \omega_{\widetilde{M}})\in \mathfrak{B}_{\widetilde{M}}$ and its image under (\ref{Isomorphism}) is $\xi(c_1(\widetilde{M}, \omega_{\widetilde{M}}))\in \mathfrak{B}_{\widehat{M}_1\cs 2 M_2}$, where the elements of the latter set are of the form (\ref{Basic Classes Contradiction}). Without loss of generality, we can assume that $\xi(c_1(\widetilde{M}, \omega_{\widetilde{M}})) = K_{\widehat{M}_1} - e^1_1 - e_1^2 - \cdots - e^2_{b_2(M_2)} - e^2_{b_2(M_2)}$. In particular,\begin{equation}\label{Rephrase Basic Class}\xi(c_1(\widetilde{M}, \omega_{\widetilde{M}})) + 2 e_1^1 = K_{\widehat{M}_1} + e_1^1 - e_1^2 - \cdots - e^2_{b_2(M_2)}\in \mathfrak{B}_{\widehat{M}_1\cs 2M_2}.\end{equation}The expression (\ref{Rephrase Basic Class}) implies that $c_1(\widetilde{M}, \omega_{\widetilde{M}}) + 2\xi^{-1}(e_1^1)\in \mathfrak{B}_{\widetilde{M}}$. To reach our contradiction, we invoke results of Taubes \cite{[Taubes2], [Taubes3]} that say that the Poincar\'e dual of $\xi^{-1}(e_1^1)$ is represented by a symplectically embedded two-sphere $S\subset (\widetilde{M}, \omega_{\widetilde{M}})$ of self-intersection $[S]^2 = -1$. This implies that $\widetilde{M} = \widetilde{M}_0\cs \overline{\mathbb{CP}^2}$ by \cite[Proposition 2.2.11, Remark 10.2.5 (a)]{[GompfStipsicz]}, contradicting the assumption of minimality of $\widetilde{M}$.
The claim on the value (\ref{SW CS}) of the Seiberg-Witten invariant is an immediate consequence of Proposition \ref{Proposition SW Connected Sum}.\end{proof}

\begin{remark}The previous modification of the argument Park yields a proof of the following result, which is a slight generalization of \cite[Theorem 6]{[Park]}. 

\begin{theorem}\label{Theorem Irreducible Diffeo Cyclic Group}Let $Z$ be a closed minimal symplectic four-manifold with $b_2^+(Z) > 1$ and let $M$ be a closed four-manifold with cyclic fundamental group $\pi_1(M) = \Z/q$. Suppose that there is an isomorphism\begin{equation}\label{Isomorphism Basic}\xi: H^2(Z)\rightarrow H^2(M)\end{equation}that preserves the cup product and that restricts to a one-to-one correspondence $\mathfrak{B}_Z\rightarrow \mathfrak{B}_M$. The four-manifold $M$ is irreducible. 
\end{theorem}

\end{remark}

\subsection{Proof of Theorem \ref{Theorem Main}}Denote by $X_n(2)\mathrel{\mathop:}= \widetilde{X}_n/\tau_n$ the closed smooth oriented four-manifold with fundamental group of order two that arises as the quotient of the closed oriented simply connected four-manifold $\widetilde{X}_n$ by its orientation-preserving fixed-point free involution $\tau_n:\widetilde{X}_n\rightarrow \widetilde{X}_n$. The classification result of Hambleton-Kreck \cite[Theorem C]{[HambletonKreck]} says that $X_n(2)$ is homeomorphic to $Y\cs \Sigma_{2, k}$ for a given closed simply connected topological four-manifold $Y$ with Kirby-Siebenmann invariant $\Ks(Y) = k = \Sigma_{2, k}$. Fix an odd positive integer $p\in \Z_{>0}$ and assemble\begin{equation}\label{Exotic Example 1}X_n(2p)\mathrel{\mathop:}= (X_n(2)\setminus \nu(\alpha_{X_n(2)}))\cup S_p,\end{equation}where $\alpha_{X_n(2)}\subset X_n(2)$ is a simple loop whose homotopy class is the generator \[\langle [\alpha_{(X_n(2)}]\rangle = \pi_1(X_n(2)) = \Z/2\] and $S_p$ is the compact four-manifold that was manufactured in Section \ref{Section Oriented Trick}. Lemma \ref{Lemma pi orientable} says that the fundamental group of the four-manifold (\ref{Exotic Example 1}) is the finite cyclic group $\pi_1(X_n(2p)) = \Z/2p$. Moreover, there is a homeomorphism between $X_n(2p)$ and $Y\cs \Sigma_{2p, k}$ by Corollary \ref{Corollary Homeomorphism Orientable} (and \cite[Theorem C]{[HambletonKreck]} as well). 


Lemma \ref{Lemma Decomposition Two Cover} says that the degree-two cover $X_n(p)$ that corresponds to the subgroup $\Z/p\subset \Z/2p$ is diffeomorphic to $X_n(p) = \widetilde{X}_n\cs \Sigma_p$. Proposition \ref{Proposition SW Procedure} says that the four-manifold (\ref{Exotic Example 1}) is irreducible. At this point, we have proven all clauses of Theorem \ref{Theorem Main} except for the third one. To conclude on the existence of infinitely many pairwise non-diffeomorphic irreducible smooth structures on $Y\cs \Sigma_{2p, k}$ and conclude the proof of the theorem, we employ Lemma \ref{Lemma Lifts} and distinguish them in terms of their double covers as done in \cite{[StipsiczSzabo2]}: the value of the Seiberg-Witten function is assumed to be $\Sw_{\widetilde{X}_n}(\pm \mathfrak{s}_n) = \pm n^2$ and there is only one pair of basic classes. Hence, $\widetilde{X}_{n_1}$ is diffeomorphic to $\widetilde{X}_{n_2}$ if and only if $n_1 = n_2$.\hfill $\square$

\section{Proof of Theorem \ref{Theorem Geo}: new irreducible examples with cyclic fundamental group.}\label{Section Proof of Theorem Geo}Theorem \ref{Theorem Main} and the following result are the main ingredients in our proof of Theorem \ref{Theorem Geo}. The case $m = 1$ is completely covered in \cite[Theorem 1.1]{[BaykurStipsiczSzabo]}.\begin{theorem}\label{Theorem Raw Material Geo}(Arabadji-Morgan \cite{[ArabadjiPorter]}, Baykur-Stipsicz-Szab\'o \cite{[BaykurStipsiczSzabo]}, Stipsicz-Szab\'o \cite{[StipsiczSzabo1], [StipsiczSzabo2]}). Let $(m, n)\in \Z_{> 0}\times \Z_{> 0}$ be positive integers that satisfy the identities\begin{itemize}\item $4 + 5 m\geq n$, \item  $4 + 5 n\geq m$, \item $m\neq 0 \mod 2$ and \item $m > 7$ if $m = n$.\end{itemize}There is a closed irreducible simply connected symplectic four-manifold $\widetilde{X}_{m, n}$ with\begin{center}$b_2^+(\widetilde{X}_{m, n}) = 2m + 1$ and $b_2^-(\widetilde{X}_{m, n}) = 2n + 1$\end{center}that admits a fixed-point free orientation-preserving involution. 
\end{theorem}Theorem \ref{Theorem Geo} follows from Theorem \ref{Theorem Raw Material Geo} and Theorem \ref{Theorem Main}.\hfill $\square$


\begin{remark}It is interesting to understand what is obtained from applying our mechanism to other constructions of exotic four-manifolds like the recent ones due to Davis-Hayden-Huang-Ruberman-Sunukjian \cite{[DHHRS]} and to Harris-Naylor-Park \cite{[HarrisNaylorPark]}, for example. 
\end{remark}

\section{Constructions of nonorientable four-manifolds and the examples of Theorem \ref{Theorem D}.}\label{Section Nonorientable}

\subsection{Cut-and-paste procedure: nonorientable setting}\label{Section nonorientable trick} There are two ways to lay down the procedure described in Section \ref{Section Oriented Trick} in the nonorientable case. Notice that the framed simple loop $\alpha_M\subset M$ used in the construction (\ref{New Manifold Oriented}) can be assumed to be orientation-preserving even if the ambient four-manifold $M$ is nonorientable and the construction described in that section results in a nonorientable four-manifold $M(p)$. We now describe the second version of this construction. 

The nonorientable compact four-manifold $N_p$ of the second step of our procedure is constructed as follows; cf. Section \ref{Section Oriented Trick}. Denote the total space of the nonorientable three-sphere bundle over the circle by $S^3\simtimes S^1$ and let $\alpha\subset S^3\simtimes S^1$ be a simple loop whose homotopy class generates the infinite cyclic group $\pi_1(S^3\simtimes S^1) = \Z$. Let $p$ be an odd integer and denote by $\alpha_p\subset S^3\simtimes S^1$ the orientation-reversing simple loop representing the class $[\alpha]^p\in H_1(S^3\simtimes S^1; \Z) = \Z$. Its tubular neighborhood is diffeomorphic to the total space of the nonorientable three-disk bundle over the circle $\nu(\alpha_p) = D^3\simtimes S^1$. 

The compact nonorientable four-manifold\begin{equation}\label{Nonorientable Piece}
N_p=(S^3 \simtimes S^1) \setminus \nu(\alpha_p)
\end{equation}with boundary $\partial N_p = S^2\simtimes S^1$ is the building block of the second step of our construction procedure when the goal is to produce nonorientable examples. 

In the following, it will be useful to understand the universal cover of $N_p$, which we now describe via the following result.

\begin{lemma}\label{lnoncover}
For any odd $p\in \Z_{> 0}$, the compact nonorientable four-manifold with boundary $N_p$ defined in (\ref{Nonorientable Piece}) is double-covered\[S_p \xrightarrow{2:1} N_p\]by the orientable four-manifold $S_p$ defined in (\ref{bb}).
\end{lemma}

\begin{proof}
    The two four-manifolds $S_p$ and $N_p$ share the same universal cover. The deck transformations of $N_p$ are generated by the map\begin{equation}\label{Definition psi}\psi: \R \times (S^3 \setminus \overset{p}{\underset{i = 1}{\sqcup}} B_i) \rightarrow \R \times (S^3 \setminus \overset{p}{\underset{i = 1}{\sqcup}} B_i)\end{equation}defined by\begin{equation}(t,x) \mapsto (t+1, g_p(x)),\end{equation}where $g_p$ is the restriction of an orientation-reversing self-diffeomorphism $S^3\rightarrow S^3$ that cyclically permutes the three-balls $\{B_1, \dots, B_p\}$. The conclusion of the lemma follows from the existence of a diffeomorphism\begin{equation}S_p \rightarrow (\R \times (S^3 \setminus \overset{p}{\underset{i = 1}{\sqcup}} B_i))/ \sim_{\psi^2}.\end{equation}\end{proof}

\begin{remark}\label{positionnonor}Thinking of $S^3$ in terms of its genus one Heegaard splitting given in (\ref{heegaardS3}), we can assume that the map $g_p$ in the proof of Lemma \ref{lnoncover} is the composition of the map $f_p:(S^3 \setminus \overset{p}{\underset{i = 1}{\sqcup}} B_i)\rightarrow (S^3 \setminus \overset{p}{\underset{i = 1}{\sqcup}} B_i)$ of Remark \ref{positionb_i's} and a reflection with respect to a two-sphere that intersects the solid torus to the left of (\ref{heegaardS3}) in $S^1 \times D^1$ and the solid torus to the right one the disjoint union of two $D^2$-fibers.\end{remark}

The following result records a key property of the four-manifold (\ref{Nonorientable Piece}) that allows for the second and third step of our procedure to be successful in the nonorientable case. This property is an extension of a foundational result due to C\'esar de S\'a \cite{[CesardeSa]} and Miller-Naylor \cite{[MillerNaylor]}, which states that any self-diffeomorphism of $S^2\simtimes S^1$ extends to a self-diffeomorphism of $D^3\simtimes S^1$.

\begin{theorem}\label{Theorem CMN}(cf. Theorem \ref{Theorem Extension}) For every odd integer $p\geq 1$ and for any self-diffeomorphism\begin{equation*}f: S^2\simtimes S^1\rightarrow S^2\simtimes S^1,\end{equation*} there is a diffeomorphism\begin{equation*}F_p: N_p\rightarrow N_p\end{equation*} such that $F_p|_{\partial} = f$. 
\end{theorem}

Theorem \ref{Theorem Extension} follows from Theorem \ref{Theorem LauPoe} and Theorem \ref{Theorem CMN}. We denote by\begin{equation}\label{nonorientable gluck}\varphi: S^2\simtimes S^1 \rightarrow S^2 \simtimes S^1\end{equation}the diffeomorphism given by $\varphi(\theta, x) = (\theta, r_{\theta}(x))$, where $r_{\theta}: S^2\rightarrow S^2$ is the rotation around a fixed axis of the two-sphere fiber about an angle $\theta$ in the circle base of $S^2\widetilde{\times} S^1$.

\begin{proof}The following argument is similar to the one used to prove Theorem \ref{Theorem LauPoe}. Kim-Raymond \cite{[KimRaymond]} showed that the diffeotopy group of $S^2 \simtimes S^1$ is isomorphic to $\Z/2 \times \Z/2$, and that it is generated by a reflection and by the diffeomorphism (\ref{nonorientable gluck}). We prove the extension result for the latter self-diffeomorphisms since the case of the reflection is immediate. Without loss of generality, we work in the context of Remark \ref{positionnonor}.  Consider the map\begin{equation}G: \R \times (S^3 \setminus \overset{p}{\underset{i = 1}{\sqcup}} B_i)\rightarrow S^3 \setminus \overset{p}{\underset{i = 1}{\sqcup}} B_i\end{equation}defined as\begin{equation}(t,x) \mapsto r_t(x),\end{equation}where $r_t: S^3\rightarrow S^3$ is a rotation of an angle of $\pi t$ that restricts to a rotation along the $D^2$-factor of the solid torus $D^2 \times S^1$ of the Heegaard decomposition of $S^3$ described in Remark \ref{positionb_i's}. Let $\psi$ be the map defined in (\ref{Definition psi}) that generates the deck transformations of $N_p$ in Lemma \ref{lnoncover}. The diffeomorphism\begin{equation}\label{Diffeomorphism T}T: \R \times (S^3 \setminus \overset{p}{\underset{i = 1}{\sqcup}} B_i) \rightarrow \R \times (S^3 \setminus \overset{p}{\underset{i = 1}{\sqcup}} B_i)\end{equation}defined by\begin{equation}(t,x) \mapsto (t, G(t,x))\end{equation}is $\psi$-equivariant and, therefore, it induces a self-diffeomorphism\begin{equation}\label{Diffeo Np}F_p: N_p \rightarrow N_p.\end{equation}

To conclude the proof, it is enough to observe that the map\begin{equation}\label{Diffeomorphism Tp}T_p:(\mathbb{R}\times (S^3 \setminus \overset{p}{\underset{i = 1}{\sqcup}} B_i)/\sim_{\psi^2} \rightarrow (\mathbb{R}\times (S^3 \setminus \overset{p}{\underset{i = 1}{\sqcup}} B_i)/\sim_{\psi^2} \cong S_p\end{equation} induced by the diffeomorphism (\ref{Diffeomorphism T}) on the double cover $S_p\rightarrow N_p$ (see (\ref{Definition psi}) and Lemma \ref{lnoncover}) restricts to a Gluck twist on the boundary $\partial S_p$. Notice that $T_p$ is the lift of $F_p$ to $S_p$. The fact that such restriction of (\ref{Diffeomorphism Tp}) is indeed a Gluck twist follows from the existence of a commutative diagram\begin{equation*}\begin{tikzcd}
(\mathbb{R}\times (S^3 \setminus \overset{p}{\underset{i = 1}{\sqcup}} B_i))/ \sim_{\phi}  \arrow[rr, "\Phi_p"] \arrow[d, "m_2"] &  & (\mathbb{R}\times (S^3 \setminus \overset{p}{\underset{i = 1}{\sqcup}} B_i)) / \sim_{\phi} \arrow[d, "m_2"] \\
(\mathbb{R}\times (S^3 \setminus \overset{p}{\underset{i = 1}{\sqcup}} B_i))/ \sim_{\psi^2} \arrow[rr, "T_p"]                      &  & (\mathbb{R}\times (S^3 \setminus \overset{p}{\underset{i = 1}{\sqcup}} B_i))/ \sim_{\psi^2},\end{tikzcd}\end{equation*}where the horizontal diffeomorphism $\Phi_p$ is defined in the proof of Theorem \ref{Theorem LauPoe}, the map $\phi$ is defined in (\ref{Phi Lemma}), and\begin{equation}m_2 \colon \mathbb{R}\times (S^3 \setminus \overset{p}{\underset{i = 1}{\sqcup}} B_i) \rightarrow \mathbb{R}\times (S^3 \setminus \overset{p}{\underset{i = 1}{\sqcup}} B_i)\end{equation}is given by $m_2(t,x)=(2t,x)$. This concludes the proof of Theorem~\ref{Theorem CMN}. For the convenience of the reader, we mention that another way to argue is to notice that the restriction of the diffeomorphism  (\ref{Diffeo Np}) to the boundary $\partial N_p = S^2\simtimes S^1$ is the composition of $p$ times the diffeomorphism (\ref{nonorientable gluck}). Such composition is isotopic to (\ref{nonorientable gluck}) given that $p$ was assumed to be odd by the result of Kim-Raymond \cite{[KimRaymond]}.\end{proof}

Much like it was done in Section \ref{Section Oriented Trick}, remove the tubular neighborhood of an orientation-reversing simple loop $\alpha_Y\subset Y$ from a closed nonorientable four-manifold $Y$ to obtain a compact nonorientable four-manifold $Y\setminus \nu(\alpha_Y)$ with boundary $\partial(Y\setminus \nu(\alpha_Y)) = S^2\simtimes S^1$. Cap it off with the compact four-manifold  (\ref{Nonorientable Piece}) to construct the closed nonorientable four-manifold\begin{equation}\label{New Manifold NonOrientable}Y(p)\mathrel{\mathop:}= (Y\setminus \nu(\alpha_Y)) \cup N_p.\end{equation}

A useful consequence of Theorem \ref{Theorem CMN} is recorded in the following corollary. 

\begin{corollary}\label{Corollary Homeomorphism Nonorientable}(cf. Theorem \ref{Theorem Homeomorphism}). Let $Y_1$ and $Y_2$ be topological four-manifolds and let $\alpha_{Y_i}\subset Y_i$ be locally flat orientation-reversing simple loops, for $i=1,2$. Suppose there is a homeomorphism
\begin{equation*}
f:Y_1\setminus \nu(\alpha_{Y_1})\rightarrow Y_2\setminus \nu(\alpha_{Y_2}).
\end{equation*}
The topological four-manifolds
\[Y_1(p)\mathrel{\mathop:}= (Y_1\setminus \nu(\alpha_{Y_1}))\cup N_p \quad \text{and} \quad Y_2(p)\mathrel{\mathop:}= (Y_2\setminus \nu(\alpha_{Y_2}))\cup N_p\] 
are homeomorphic.
\end{corollary}

The simplicity of the construction makes it straight-forward to compute several basic invariants of the new four-manifold (\ref{New Manifold NonOrientable}) as we now sample.

\begin{lemma}\label{Lemma pi nonorientable}The closed nonorientable four-manifold $Y(p)$ defined in (\ref{New Manifold NonOrientable}) has the following basic topological invariants for every $p\in \Z_{> 0}$.\begin{enumerate}\item The Euler characteristic is $\chi(Y(p)) = \chi(Y)$.
\item The second Stiefel-Whitney class is $w_2(Y(p)) = w_2(Y)$ and the relation $w_1(Y(p))^2 + w_2(Y(p)) = w_1(Y)^2 + w_2(Y)$ holds.
\item If the homotopy class of the simple loop $\alpha_Y\subset Y$ generates the cyclic fundamental group $\pi_1(Y) = \Z/q$, then the fundamental group of $Y(p)$ is $\pi_1(Y(p)) = \langle t: t^{q\cdot p} = 1\rangle$.
\end{enumerate} 
\end{lemma}

Much like the in the case of Lemma \ref{Lemma pi orientable}, a proof of Lemma \ref{Lemma pi nonorientable} is obtained through straight-forward (co)homological computations, cobordism arguments and the Seifert-van Kampen theorem.

\subsection{Nonorientable raw material and its properties}\label{Section NonOrientable Background}With the purpose of streamlining the exposition, we encapsule the background results we build on to prove Theorem \ref{Theorem D} in the next triad of results. We first recall the building block that serves as raw material of the procedure we described in the introduction, along with several of its basic properties. 

\begin{theorem}\label{Theorem NonOrientable 1}\

$\bullet$ Cappell-Shaneson \cite{[CappellShaneson]}. There is a four-manifold $Q$ that is homeomorphic, but not diffeomorphic to $\RP^4$. These two four-manifolds remain non-diffeomorphic under connect sums with arbitrarily many copies of $S^2\times S^2$. 

$\bullet$ Gompf \cite{[Gompf]}. The universal cover $\widetilde{Q}$ is diffeomorphic to $S^4$.

$\bullet$ Akbulut \cite[Theorem 2]{[Akbulut]}. There is a diffeomorphism\begin{equation*}Q\cs \CP^2\rightarrow \RP ^4\cs \CP^2.\end{equation*} 

\end{theorem}

We now add a few words of the foundational results summarized in the previous theorem in order to put them into the context of our proof of Theorem \ref{Theorem D}. Recall that the tubular neighborhood of\begin{equation*}\RP^2 = \{[x:y:z:0:0] | x^2 + y^2 + z^2 = 1\}\subset S^4/{\sim}  = \RP^4\end{equation*}is the two-disk bundle over the real projective plane\begin{equation*}D^2\simtimes \RP^2 = (D^2\times S^2)/{\sim},\end{equation*}where $(d, x)\sim (-d, \ant(x))$ for every $d\in D^2$ and $x\in S^2$ and $\ant: S^2\rightarrow S^2$ is the antipodal involution. The analysis done by Akbulut \cite[Section 1]{[Akbulut]} on the four-manifold $Q$ built by Cappell-Shaneson \cite{[CappellShaneson]} reveals that the compact four-manifold $Z\mathrel{\mathop:}= Q\setminus \nu(\alpha_Q)$ (where $\alpha_Q\subset Q$ is a simple loop that generates $\pi_1(Q) = \Z/2$) is homeomorphic, but not diffeomorphic to $D^2\simtimes \RP^2$. The second clause of Theorem \ref{Theorem NonOrientable 1} encapsules work of Gompf, where he showed that the orientation two-cover of $Z$ is diffeomorphic to $D^2\times S^2$. The third clause of Theorem \ref{Theorem NonOrientable 1} follows from the property $Z\cs \CP^2 = (D^2\simtimes \RP^2)\cs \CP^2$ \cite[Section 3]{[Akbulut]}.

\subsection{A spectral invariant as a diffeomorphism obstruction}The $\eta$-invariant will be used to distinguish the diffeomorphism type of our nonorientable four-manifolds. We refer the reader to \cite{[Stolz]} for further details on this spectral invariant and to \cite{[HambletonKreckTeichner], [KirbyTaylor], [Stolz]} for background on $\Pin^+$-structures. We say that two four-manifolds $Y$ and $Y'$ are stably diffeomorphic if there is a positive integer $k\in \Z_{> 0}$ such that there is a diffeomorphism between the connected sums\begin{center}$Y\cs (k - 1)(S^2\times S^2)$ and $Y'\cs (k - 1)(S^2\times S^2)$,\end{center}where we set $0(S^2\times S^2) = S^4$.

\begin{theorem}[Stolz \cite{[Stolz]}]\label{Theorem NonOrientable 2}\

$\bullet$ If the $\eta$-invariant of two closed nonorientable $\Pin^+$- four-manifolds $(Y, \phi_Y)$ and $(Y', \phi_{Y'})$ with $H^1(Y; \Z/2) = \Z/2 = H^1(Y'; \Z/2)$ satisfies $\eta(Y, \phi_Y) \neq \eta(Y', \phi_{Y'}) \mod 2\Z$, then $Y_1$ is not stably diffeomorphic to $Y_2$. 

$\bullet$ There are $\Pin^+$-structures $(\RP^4, \pm \phi_{\RP^4})$ and $(Q, \pm \phi_Q)$ such that\begin{equation*}\eta(\RP^4, \pm \phi_{\RP^4}) = \pm \frac{1}{8} \mod 2\Z\end{equation*}and\begin{equation*}\eta(Q, \pm \phi_Q) = \pm \frac{7}{8} \mod 2\Z.\end{equation*}In particular, $Q$ is not stably diffeomorphic to $\RP^4$. 

\end{theorem}

The last claim  was proven by Cappell-Shaneson in \cite{[CappellShaneson]}.

\subsection{$\Pin^+$-cobordism classes and a computational tool}\label{Section Computations Eta}We now elaborate on the role that the $\eta$-invariant plays in telling two given smooth structures apart. As we mentioned in the introduction, we distinguish the diffeomorphism classes by looking at the $\Pin^+$-structures on the tangent bundle of our four-manifolds through the $\eta$-invariant. 

\begin{theorem}\label{Theorem NonOrientable 3}\

$\bullet$ Kirby-Taylor \cite[Theorem 5.2]{[KirbyTaylor]}. There is a group isomorphism \[\Omega^{\Pin^+}_4\rightarrow \Z/16\] and the class of the real projective four-space\begin{equation*}[(\RP^4, \pm \phi_{\RP^4})] = \pm 1\in \Omega^{\Pin^+}_4\cong \Z/16.\end{equation*}is a generator.

$\bullet$ Stolz \cite{[Stolz]}. The $\eta$-invariant$\mod 2\Z$ is a complete $\Pin^+$-cobordism invariant and it is additive under disjoint unions. 

$\bullet$ Stolz \cite{[Stolz]}.\begin{equation*}[(Q, \pm \phi_{Q})] = \pm 9\in \Omega^{\Pin^+}_4\cong \Z/16.\end{equation*}
\end{theorem}

We finish this section with the main computational tool to distinguish the four-manifolds of our procedure in the presence of a $\Pin^+$-structure. 

\begin{proposition}\label{Proposition Pin Cobordisms}Let $Y$ be a closed four-manifold with first cohomology group $H^1(Y; \Z/2) = \Z/2$ and that admits a $\Pin^+$-structure $(Y, \phi_Y)$. Let $\alpha_Y\subset Y$ be an orientation-reversing simple loop such that $[\alpha_Y]\neq 0\in H_1(Y; \Z)$ and build\begin{equation}\label{Construct Cob}Y(p)\mathrel{\mathop:}= (Y\setminus \nu(\alpha_Y))\cup N_p.\end{equation}

There are $\Pin^+$-structures\begin{center}$(S^3\simtimes S^1, \phi)$ and $(Y(p), \phi_{Y(p)})$\end{center}along with a $\Pin^+$-cobordism\begin{equation}\label{Pin Cobordism N}((V, \phi_V); (Y, \phi_Y)\sqcup(S^3\simtimes S^1, \phi), (Y(p), \phi_{Y(p)})).\end{equation}In particular, the identity\begin{equation}\label{Identity Cobordism}[(Y, \phi_Y)] = [(Y(p), \phi_{Y(p)})]\in \Omega^{\Pin^+}_4\end{equation}holds.  
\end{proposition}

A similar statement to Proposition \ref{Proposition Pin Cobordisms} continues to hold if the simple loop $\alpha_Y\subset Y$ is assumed to be orientation-preserving.

\begin{proof}We first assemble the $\Pin^+$-structure $(Y(p), \phi_{Y(p)})$ by gluing together $\Pin^+$-structures $(Y\setminus \nu(\alpha_Y), \phi_Y^0)$ and $(N_p, \phi_{N_p})$ that induce the same $\Pin^+$-structure on their common boundary; see (\ref{Construct Cob}). The $\Pin^+$-structure $(Y, \phi_Y)$ induces a $\Pin^+$-structure $(Y\setminus \nu(\alpha_Y), \phi_Y^0)$ as a codimension zero submanifold \cite[Construction 1.5]{[KirbyTaylor]}, which restricts to a  $\Pin^+$-structure\begin{equation}\label{Choice Pin 1}(\partial (Y\setminus \nu(\alpha_Y)), \phi^0_Y|_{\partial}) = (S^2\simtimes S^1, \phi_\partial).\end{equation} Since\begin{equation*}H^1(Y; \Z/2) = \Z/2 = H^1(Y\setminus \nu(\alpha_Y); \Z/2) = H^1(\partial (Y\setminus \nu(\alpha_Y)); \Z/2),\end{equation*}there are two choices of such structures \cite[Corollary 6.4]{[Stolz]} and there is a one-to-one correspondence between the initial $\Pin^+$-structure $(Y, \phi_Y)$ and the $\Pin^+$-structure (\ref{Choice Pin 1}) (see \cite[Proposition 1 Chapter IV]{[Kirby]}). We now argue that there is a $\Pin^+$-structure $(N_p, \phi_{N_p})$ that induces the $\Pin^+$-structure (\ref{Choice Pin 1}) on $\partial N_p = S^2\simtimes S^1$. Since $w_2(S^3\simtimes S^1) = 0$ and $H^1(S^3\simtimes S^1; \Z/2) = \Z/2$, there are two $\Pin^+$-structures on $S^3\simtimes S^1$ \cite[Corollary 6.4]{[Stolz]}. The hypothesis $p\neq 0 \mod 2$ allows us to pick the $\Pin^+$-structure $(S^3\simtimes S^1, \phi)$ that induces $(N_p, \phi_{N_p})$ such that the induced $\Pin^+$-structure $(\partial N_p, \phi_{N_p}|_\partial)$ matches (\ref{Choice Pin 1}). We conclude that there is a $\Pin^+$-structure $(Y(p), \phi_{Y(p)})$. 
Next, we put together the compact five-manifold $V$ by gluing the product cobordisms\begin{center}$Y\times [0, 1]$ and $(S^3\simtimes S^1)\times [0, 1]$\end{center}along the upper lids $Y\times \{1\}$ and $(S^3\simtimes S^1)\times \{1\}$ according to the construction of $Y(p)$ in (\ref{Construct Cob}). The $\Pin^+$-structure $(V, \phi_V)$ is assembled analogously to the assembly of $(Y, \phi_Y)$ that was described before. Given that\begin{equation*}[(S^3\simtimes S^1, \pm \phi)] = 0\in \Omega^{\Pin^+}_4,\end{equation*}the identity (\ref{Identity Cobordism}) follows from the existence of the cobordism (\ref{Pin Cobordism N}).\end{proof}

\subsection{Proof of Theorem \ref{Theorem D}}\label{Section Proof Theorem E}Fix an odd $p\in \Z_{\geq 0}$ and consider the compact four-manifold $N_p$ that was defined in (\ref{Nonorientable Piece}). Assemble the nonorientable four-manifolds\begin{equation}\label{Manifold A}A(2p)\mathrel{\mathop:}= (\RP^4\setminus \nu(\alpha_{\RP^4}))\cup N_p = (D^2\simtimes \RP^2)\cup N_p\end{equation}and\begin{equation}\label{Manifold B}B(2p)\mathrel{\mathop:}= (Q\setminus \nu(\alpha_Q))\cup N_p = Z\cup N_p,\end{equation}where $Q$ is Cappell-Shaneson's exotic copy of the real projective four-space \cite{[CappellShaneson]} discussed in Theorem \ref{Theorem NonOrientable 1}. Lemma \ref{Lemma pi nonorientable} implies that the four-manifolds (\ref{Manifold A}) and (\ref{Manifold B}) are $\Q$-homology real projective four-spaces whose fundamental group is a finite cyclic group of order $2p$. This establishes the claim of Item (1) of Theorem \ref{Theorem D}. Moreover, the four-manifolds $A(2p)$ and $B(2p)$ are homeomorphic by Corollary \ref{Corollary Homeomorphism Nonorientable} since $Q$ is homeomorphic to $\RP^4$.  This settles the claim of Item (2). 
Since there are $\Pin^+$-structures $(\RP^4, \pm \phi_{\RP^4})$ and $(Q, \pm \phi_Q)$, Proposition \ref{Proposition Pin Cobordisms} implies that the four-manifolds (\ref{Manifold A}) and (\ref{Manifold B}) admit a $\Pin^+$-structure. There are actually two such structures. Indeed, the first cohomology groups with $\Z/2$-coefficients of the four-manifolds (\ref{Manifold A}) and (\ref{Manifold B}) are\begin{center}$H^1(A(2p); \Z/2) = \Z/2 = H^1(B(2p); \Z/2)$\end{center}and it follows that there are two $\Pin^+$-structures $(A(2p), \pm \phi_{A(2p)})$ and two $\Pin^+$-structures $(B(2p), \pm \phi_{B(2p)})$ \cite[Corollary 6.4]{[Stolz]}, . Moreover, the $\Pin^+$-structures on the same underlying four-manifold are mutual inverses in the fourth $\Pin^+$-cobordism group $\Omega^{\Pin^+}_4$ and the values\begin{equation}
\eta(A(2p), \pm \phi_{A(2p)}) = \pm \frac{1}{8} \mod 2\Z
\end{equation}and\begin{equation}
\eta(B(2p), \pm \phi_{B(2p)}) = \pm \frac{7}{8} \mod 2\Z
\end{equation}
hold by Proposition \ref{Proposition Pin Cobordisms}, Theorem \ref{Theorem NonOrientable 2} and Theorem \ref{Theorem NonOrientable 3}. Using Theorem \ref{Theorem NonOrientable 2}, we conclude that the connected sums\begin{center}$A(2p)\cs (k - 1)(S^2\times S^2)$ and $B(2p)\cs (k - 1)(S^2\times S^2)$\end{center} belong to different classes in $\Omega^{\Pin^+}_4\cong \Z/16$ for any $k\in \Z_{> 0}$ and that the four-manifolds $A(2p)$ and $B(2p)$ are not stably-diffeomorphic. This settles the claim of Item (5). Item (6) follows from the third clause of Theorem \ref{Theorem NonOrientable 1} and the corresponding covers of the decompositions (\ref{Manifold A}) and (\ref{Manifold B}). 

Finally, to construct the closed non-smoothable topological four-manifold $\ast A(2p)$, we substitute $\RP^4$ with the closed non-smootable four-manifold $\ast \RP^4$ that is homotopy equivalent to the real projective four-space in (\ref{Manifold A}). The latter four-manifold exists by a computation using the surgery exact sequence \cite[Chapter 14]{[Wall]} and an explicit construction of it was given by Ruberman in \cite{[Ruberman]}. The existence and uniqueness results of a topological tubular neighborhood $\nu(\alpha_{\ast \RP^4})$ for a locally flat simple loop $\alpha_{\ast \RP^4}\subset \ast \RP^4$ are found in \cite[Section 5.1]{[FNOP]}. Choose a loop that generates the fundamental group and assemble\begin{equation}\label{Manifold ANS}\ast A(2p)\mathrel{\mathop:}= (\ast\RP^4\setminus \nu(\alpha_{\ast\RP^4}))\cup N_p.\end{equation}The existence of a homotopy equivalence $\ast \RP^4\rightarrow \RP^4$ implies that there is a homotopy equivalence $\ast A(2p)\rightarrow A(2p)$. From (\ref{Manifold ANS}) and the related cobordism, we have that the Kirby-Siebenmann invariant is $\Ks(\ast A(2p)) = \Ks(\ast \RP^4) = 1\in \Z/2$  \cite[Theorem 9.2]{[FNOP]} (see \cite[Section 9.2]{[FNOP]} for a description of the construction of this invariant). Thus, the four-manifold $\ast A(2p)$ is non-smoothable for any odd $p\in \Z_{>0}$. The homeomorphism class of its universal cover is pinned down by applying Freedman's theorem \cite[Theorem 1.2.27]{[GompfStipsicz]}.\hfill $\square$

\end{document}